\documentclass[12pt]{article}
\usepackage[]{amsmath,amssymb}
\usepackage{latexsym}
\usepackage{cite}

\newtheorem{definition}{Definition}[section]
\newtheorem{theorem}[definition]{Theorem}
\newtheorem{lemma}[definition]{Lemma}

\newtheorem{example}[definition]{Example}

\newtheorem{note}[definition]{Note}
\newtheorem{assumption}[definition]{Assumption}
\newtheorem{proposition}[definition]{Proposition}

\def\F{\mathbb F}
\def\K{\mathbb F}

\begin{document}

\title{\bf The Drinfel'd polynomial  \\
of a tridiagonal pair
}
\author{
Tatsuro Ito
 and
Paul Terwilliger\footnote{This author gratefully acknowledges 
support from the FY2007 JSPS Invitation Fellowship Program
for Reseach in Japan (Long-Term), grant L-07512.}
}
\date{}

\maketitle

\centerline{\large \bf Dedicated to Dijen Ray-Chaudhuri on 
his 75th birthday}

\begin{abstract}
Let $\K$ denote a field and let $V$ denote a vector space
over $\K$ with finite positive dimension.
We consider a pair of linear transformations $A:V \to V$
and $A^*:V \to V$ that satisfy the following conditions:
(i)
each of $A,A^*$ is diagonalizable;
(ii)
there exists an ordering $\lbrace V_i\rbrace_{i=0}^d$ of the eigenspaces of
$A$ such that
$A^* V_i \subseteq V_{i-1} + V_{i} + V_{i+1}$ for $0 \leq i \leq d$,
where $V_{-1}=0$ and $V_{d+1}=0$;
(iii)
there exists an ordering $\lbrace V^*_i\rbrace_{i=0}^\delta$
of the eigenspaces of $A^*$ such that
$A V^*_i \subseteq V^*_{i-1} + V^*_{i} + V^*_{i+1}$ for
 $0 \leq i \leq \delta$,
where $V^*_{-1}=0$ and $V^*_{\delta+1}=0$;
(iv)
there is no subspace $W$ of $V$ such that
$AW \subseteq W$, $A^* W \subseteq W$, $W \neq 0$, $W \neq V$.
We call such a pair a {\it tridiagonal pair} on $V$.
It is known that $d=\delta$ and for $0 \leq i \leq d$
the dimensions of $V_i$, $V_{d-i}$, $V^*_i$, $V^*_{d-i}$ coincide.
The pair $A,A^*$ is called {\it sharp} whenever $\dim V_0=1$.
It is known that if
$\K$ is algebraically closed then $A,A^*$ is sharp.
Assuming $A,A^*$ is sharp, we use the
data 
$\Phi=(A; 
\lbrace V_i\rbrace_{i=0}^d;
A^*; 
\lbrace V^*_i\rbrace_{i=0}^d)$
 to define
a polynomial
$P$ in one variable and degree at most $d$.
We show that $P$
remains invariant if $\Phi$ is replaced by
$(A;\lbrace V_{d-i}\rbrace_{i=0}^d;
A^*; 
\lbrace V^*_i\rbrace_{i=0}^d)$
or 
$(A;\lbrace V_i\rbrace_{i=0}^d;
A^*; 
\lbrace V^*_{d-i}\rbrace_{i=0}^d)$
or 
$(A^*; 
\lbrace V^*_i\rbrace_{i=0}^d;
A; 
\lbrace V_i\rbrace_{i=0}^d)$.
We call $P$ the {\it Drinfel'd polynomial}
of $A,A^*$.  
We explain how $P$ is related to the
classical Drinfel'd polynomial from the
theory of Lie algebras and quantum groups.
We expect that the roots of $P$ will 
be useful in a future classification of the sharp tridiagonal pairs.
We compute the roots of $P$ for the case in
which $V_i$ and $V^*_i$ have dimension
1 for $0 \leq i \leq d$.

\bigskip
\noindent
{\bf Keywords}. 
Tridiagonal pair, Leonard pair, $q$-Racah polynomial.
 \hfil\break
\noindent {\bf 2000 Mathematics Subject Classification}. 
Primary: 15A21. Secondary: 
05E30, 05E35,
33D45, 17B37.
 \end{abstract}

\section{Tridiagonal pairs}

\noindent 
Throughout this paper $\K$ denotes a field
 and
$V$ denotes a vector space over $\K$ with finite
positive dimension.

\medskip
\noindent 
We begin by recalling the notion of a tridiagonal pair. 
We will use the following terms.
For a linear transformation $A:V\to V$ 
and a  subspace $W \subseteq V$,
we call $W$ an
 {\it eigenspace} of $A$ whenever 
 $W\not=0$ and there exists $\theta \in \K$ such that 
$W=\lbrace v \in V \;\vert \;Av = \theta v\rbrace$;
in this case $\theta$ is the {\it eigenvalue} of
$A$ associated with $W$.
We say that $A$ is {\it diagonalizable} whenever
$V$ is spanned by the eigenspaces of $A$.

\begin{definition}  
{\rm \cite[Definition~1.1]{TD00}}
\label{def:tdp}
\rm
By a {\it tridiagonal pair} 
on $V$
we mean an ordered pair of linear transformations
$A:V\to V$ and 
$A^*:V\to V$ 
that satisfy the following four conditions.
\begin{enumerate}
\item Each of $A,A^*$ is diagonalizable.
\item There exists an ordering $\lbrace V_i\rbrace_{i=0}^d$ of the  
eigenspaces of $A$ such that 
\begin{equation}
A^* V_i \subseteq V_{i-1} + V_i+ V_{i+1} \qquad \qquad (0 \leq i \leq d),
\label{eq:t1}
\end{equation}
where $V_{-1} = 0$ and $V_{d+1}= 0$.
\item There exists an ordering $\lbrace V^*_i\rbrace_{i=0}^{\delta}$ of
the  
eigenspaces of $A^*$ such that 
\begin{equation}
A V^*_i \subseteq V^*_{i-1} + V^*_i+ V^*_{i+1} 
\qquad \qquad (0 \leq i \leq \delta),
\label{eq:t2}
\end{equation}
where $V^*_{-1} = 0$ and $V^*_{\delta+1}= 0$.
\item There does not exist a subspace $W$ of $V$ such  that $AW\subseteq W$,
$A^*W\subseteq W$, $W\not=0$, $W\not=V$.
\end{enumerate}
We say the pair $A,A^*$ is {\it over $\K$}.
\end{definition}

\begin{note}
\rm
According to a common notational convention $A^*$ denotes 
the conjugate-transpose of $A$. We are not using this convention.
In a tridiagonal pair $A,A^*$ the linear transformations $A$ and $A^*$
are arbitrary subject to (i)--(iv) above.
\end{note}

\noindent We refer the reader to 
\cite{hasan,
hasan2,
shape,
tdanduq,
NN,
N:refine,
nomsplit,
qSerre,
madrid} and the references therein 
for background on tridiagonal pairs. 

\medskip
\noindent In order to motivate our results we
recall a few facts about tridiagonal pairs.
Let $A,A^*$ denote a tridiagonal pair
on $V$, as in Definition 
\ref{def:tdp}. By
\cite[Lemma 4.5]{TD00}
the integers $d$ and $\delta$ from
(ii), (iii) are equal; we call this
common value the {\it diameter} of the
pair.
An ordering of the eigenspaces of $A$ (resp. $A^*$)
is said to be {\it standard} whenever it satisfies 
(\ref{eq:t1})
 (resp. (\ref{eq:t2})). 
We comment on the uniqueness of the standard ordering.
Let $\{V_i\}_{i=0}^d$ denote a standard ordering of the eigenspaces of $A$.
Then the ordering $\{V_{d-i}\}_{i=0}^d$ is also standard and no
 further ordering
is standard.
A similar result holds for the eigenspaces of $A^*$.
Let $\{V_i\}_{i=0}^d$ (resp.
$\{V^*_i\}_{i=0}^d$)
denote a standard ordering of the eigenspaces of
 $A$ (resp. $A^*$).
By \cite[Corollary 5.7]{TD00}, 
for $0 \leq i \leq d$ the spaces $V_i$, $V^*_i$
have the same dimension; we denote
this common dimension by $\rho_i$. 
By \cite[Corollaries 5.7, 6.6]{TD00}
the sequence $\{\rho_i\}_{i=0}^d$ is symmetric and unimodal;
that is $\rho_i=\rho_{d-i}$ for $0 \leq i \leq d$ and
$\rho_{i-1} \leq \rho_i$ for $1 \leq i \leq d/2$.
We call the sequence $\{\rho_i\}_{i=0}^d$ the {\em shape}
of $A,A^*$. We say $A,A^*$ is {\it sharp} whenever
$\rho_0=1$.
By 
\cite[Theorem~1.3]{nomstructure},
if $\K$ is algebraically closed then $A,A^*$ is 
sharp.
It is an open problem to  classify the sharp tridiagonal pairs
up to isomorphism
\cite{nommu}.
By a {\em Leonard pair} we mean a tridiagonal pair with shape
$(1,1,\ldots,1)$ \cite[Definition 1.1]{LS99}.
In
\cite{LS99,
 TLT:array}
the Leonard pairs are classified  
up to isomorphism.
This classification yields a correspondence between the Leonard pairs and a
family of orthogonal polynomials consisting of the $q$-Racah polynomials
and their relatives 
\cite{AWil,
 TLT:array,
qrac}.
This family coincides with the terminating branch of the Askey scheme
\cite{KoeSwa}.
See 
\cite{NT:formula,NT:det,NT:mu,
NT:span,NT:switch,
   conform,
    lsint,
TLT:split,
aw}
for more information about Leonard pairs.

\medskip
\noindent
We now describe our main results.
Let $A,A^*$ denote a sharp tridiagonal pair on $V$.
Let $\{V_i\}_{i=0}^d$ (resp.
$\{V^*_i\}_{i=0}^d$)
denote a standard ordering of the eigenspaces of
 $A$ (resp. $A^*$).
Using the data
$(A; 
\lbrace V_i\rbrace_{i=0}^d;
A^*; 
\lbrace V^*_i\rbrace_{i=0}^d)$
and the following construction,
we define a polynomial 
$P$ in one variable $\lambda$.
For $0 \leq i \leq d$ let $\theta_i$
(resp. $\theta^*_i$)
denote the eigenvalue of $A$ (resp. $A^*$) associated with
 $V_i$
(resp. $V^*_i$).
By \cite[Theorem~4.6]{TD00},
for $0 \leq i \leq d$ the space $V^*_0$ is invariant under
\begin{eqnarray*}
 (A^*-\theta^*_1I)(A^*-\theta^*_2I)\cdots(A^*-\theta^*_iI)
 (A-\theta_{i-1}I)\cdots(A-\theta_1I)(A-\theta_0I);
\end{eqnarray*}
let $\zeta_i$ denote the corresponding eigenvalue.
By \cite[Theorem~11.1]{TD00} there exists  $\beta \in \K$
such that $\theta_{i-1}-\beta \theta_i + \theta_{i+1}$
 and $\theta^*_{i-1}-\beta \theta^*_i + \theta^*_{i+1}$
are independent of $i$  for $1 \leq i \leq d-1$.
For the moment assume $\beta \not=\pm 2$ and put
$q^2+q^{-2}=\beta$.
Then $q^{2i}\not=1$ for $1 \leq i \leq d$,
by
Note 
\ref{lem:l3I} below.
Define
\begin{eqnarray*}
P = \sum_{i=0}^d \zeta_i p_{i+1}p_{i+2}\cdots p_d,
\end{eqnarray*}
where 
\begin{eqnarray*}
p_i = (\theta_0 \theta^*_d+\theta_d \theta^*_0-\lambda)
\frac{(q^i-q^{-i})^2}{(q^d-q^{-d})^2}
+
(\theta_0-\theta_i) 
(\theta^*_0-\theta^*_i)
\end{eqnarray*}
for $1 \leq i \leq d$.
For $\beta=\pm 2$ our definition of $P$ is slightly different;
see Definition 
 \ref{def:aap}
and 
Definition
\ref{def:3p} for the details.
For all $\beta$ we show that $P$ remains invariant
if the data  
$(A;\lbrace V_{i}\rbrace_{i=0}^d;
A^*; 
\lbrace V^*_i\rbrace_{i=0}^d)$
is replaced by
$(A;\lbrace V_{d-i}\rbrace_{i=0}^d;
A^*; 
\lbrace V^*_i\rbrace_{i=0}^d)$
or 
$(A;\lbrace V_i\rbrace_{i=0}^d;
A^*; 
\lbrace V^*_{d-i}\rbrace_{i=0}^d)$
or 
$(A^*; 
\lbrace V^*_i\rbrace_{i=0}^d;
A; 
\lbrace V_i\rbrace_{i=0}^d)$.
We call $P$ the {\it Drinfel'd polynomial}
of $A,A^*$. 
We explain how $P$ is related to the
classical Drinfel'd polynomial from the
theory of Lie algebras and quantum groups.
We expect that the roots of $P$ will be useful
 in a future classification of the sharp tridiagonal pairs.
We compute the roots of $P$ for the case in
which $A,A^*$ is a Leonard pair.

\section{Tridiagonal systems}

\indent
When working with a tridiagonal pair, it is often convenient to consider
a closely related object called a tridiagonal system.
To define a tridiagonal system, we recall a few concepts from linear
algebra.
Let ${\rm End}(V)$ denote the $\K$-algebra of all linear
transformations from $V$ to $V$.
Let $A$ denote a diagonalizable element of $\mbox{\rm End}(V)$.
Let $\{V_i\}_{i=0}^d$ denote an ordering of the eigenspaces of $A$
and let $\{\theta_i\}_{i=0}^d$ denote the corresponding ordering of
the eigenvalues of $A$.
For $0 \leq i \leq d$ define $E_i \in 
\mbox{\rm End}(V)$ 
such that $(E_i-I)V_i=0$ and $E_iV_j=0$ for $j \neq i$ $(0 \leq j \leq d)$.
Here $I$ denotes the identity of $\mbox{\rm End}(V)$.
We call $E_i$ the {\em primitive idempotent} of $A$ corresponding to $V_i$
(or $\theta_i$).
Observe that
(i) $\sum_{i=0}^d E_i = I$;
(ii) $E_iE_j=\delta_{i,j}E_i$ $(0 \leq i,j \leq d)$;
(iii) $V_i=E_iV$ $(0 \leq i \leq d)$;
(iv) $A=\sum_{i=0}^d \theta_i E_i$.
Moreover
\begin{equation}         \label{eq:defEi}
  E_i=\prod_{\stackrel{0 \leq j \leq d}{j \neq i}}
          \frac{A-\theta_jI}{\theta_i-\theta_j}.
\end{equation}
We note that each of $\{E_i\}_{i=0}^d$,
$\{A^i\}_{i=0}^d$ 
is a basis for the $\K$-subalgebra
of $\mbox{\rm End}(V)$ generated by $A$.
Now let $A,A^*$ denote a tridiagonal pair on $V$.
An ordering of the primitive idempotents 
or eigenvalues of $A$ (resp. $A^*$)
is said to be {\em standard} whenever
the corresponding ordering of the eigenspaces of $A$ (resp. $A^*$)
is standard.

\begin{definition}
{\rm \cite[Definition~2.1]{TD00}}
  \label{def:TDsystem} 
\rm
By a {\it tridiagonal system} on $V$ we mean a sequence
\[
 \Phi=(A;\{E_i\}_{i=0}^d;A^*;\{E^*_i\}_{i=0}^d)
\]
that satisfies (i)--(iii) below.
\begin{itemize}
\item[(i)]
$A,A^*$ is a tridiagonal pair on $V$.
\item[(ii)]
$\{E_i\}_{i=0}^d$ is a standard ordering
of the primitive idempotents of $A$.
\item[(iii)]
$\{E^*_i\}_{i=0}^d$ is a standard ordering
of the primitive idempotents of $A^*$.
\end{itemize}
We say $\Phi$ is {\em over} $\K$.
We call $d$ the {\it diameter} of $\Phi$.
\end{definition}

\medskip

\noindent 
The following result is immediate from lines (\ref{eq:t1}),
 (\ref{eq:t2})
and Definition \ref{def:TDsystem}.
\medskip

\begin{lemma}  
{\rm \cite[Lemma~2.5]{nomtowards}}
  \label{lem:trid}     \samepage
Let $(A;\{E_i\}_{i=0}^d;A^*;\{E^*_i\}_{i=0}^d)$ denote a tridiagonal system.
Then for $0 \leq i,j,k \leq d$ the following {\rm (i)}, {\rm (ii)} hold.
\begin{itemize}
\item[\rm (i)]
$E^*_i A^k E^*_j =0\;\;$ if $k<|i-j|$.
\item[\rm (ii)]
$E_i A^{*k} E_j = 0\;\;$ if  $k<|i-j|$.
\end{itemize}
\end{lemma}

\begin{definition}        \label{def0}
\rm
Let 
$\Phi=(A;\{E_i\}_{i=0}^d;A^*$; $\{E^*_i\}_{i=0}^d)$ 
denote a tridiagonal system.
By the {\it associated tridiagonal pair} we mean
the pair $A,A^*$.
By the {\it shape} of $\Phi$ we mean the shape
of $A,A^*$.
We say $\Phi$ is
{\it sharp} whenever 
$A,A^*$ is sharp.
We call $\Phi$ a {\it Leonard system} whenever
$A,A^*$ is a Leonard pair.
\end{definition}

\section{The $D_4$ action}
\indent
Let $\Phi=(A; \{E_i\}_{i=0}^d; A^*; \{E^*_i\}_{i=0}^d)$
denote a tridiagonal system on $V$.
Then each of the following is a tridiagonal system on $V$:
\begin{eqnarray*}
\Phi^{*}  &:=&
       (A^*; \{E^*_i\}_{i=0}^d; A; \{E_i\}_{i=0}^d), \\
\Phi^{\downarrow} &:=&
       (A; \{E_i\}_{i=0}^d; A^*; \{E^*_{d-i}\}_{i=0}^d), \\
\Phi^{\Downarrow} &:=&
       (A; \{E_{d-i}\}_{i=0}^d; A^*; \{E^*_{i}\}_{i=0}^d).
\end{eqnarray*}
Viewing $*, \downarrow, \Downarrow$
as permutations on the set of all tridiagonal systems,
\begin{eqnarray}
&&\qquad \qquad \qquad  *^2 \;=\;  
\downarrow^2\;= \;
\Downarrow^2 \;=\;1,
\qquad \quad 
\label{eq:8rel1}
\\
&&\Downarrow *\; 
=\;
* \downarrow,\qquad \qquad   
\downarrow *\; 
=\;
* \Downarrow,\qquad \qquad   
\downarrow \Downarrow \; = \;
\Downarrow \downarrow.
\qquad \quad 
\label{eq:8rel2}
\end{eqnarray}
The group generated by symbols $*$, $\downarrow$, $\Downarrow$ subject
to the relations (\ref{eq:8rel1}), (\ref{eq:8rel2}) is the
dihedral group $D_4$. We recall that $D_4$ is the group of symmetries of a
square, and has $8$ elements.
Apparently $*$, $\downarrow$, $\Downarrow$ induce an action of $D_4$
on the set of all tridiagonal systems.
Two tridiagonal systems will be called {\em relatives} whenever they are
in the same orbit of this $D_4$ action. 
The relatives of $\Phi$ are as follows:
\medskip
\noindent
\begin{center}
\begin{tabular}{c|c}
name  &  relative \\
\hline
$\Phi$ & 
       $(A; \{E_i\}_{i=0}^d; A^*;  \{E^*_i\}_{i=0}^d)$ \\ 
$\Phi^{\downarrow}$ &
       $(A; \{E_i\}_{i=0}^d; A^*;  \{E^*_{d-i}\}_{i=0}^d)$ \\ 
$\Phi^{\Downarrow}$ &
       $(A; \{E_{d-i}\}_{i=0}^d; A^*;  \{E^*_i\}_{i=0}^d)$ \\ 
$\Phi^{\downarrow \Downarrow}$ &
       $(A; \{E_{d-i}\}_{i=0}^d; A^*;  \{E^*_{d-i}\}_{i=0}^d)$ \\ 
$\Phi^{*}$  & 
       $(A^*; \{E^*_i\}_{i=0}^d; A;  \{E_i\}_{i=0}^d)$ \\ 
$\Phi^{\downarrow *}$ &
       $(A^*; \{E^*_{d-i}\}_{i=0}^d; A;  \{E_i\}_{i=0}^d)$ \\ 
$\Phi^{\Downarrow *}$ &
       $(A^*; \{E^*_i\}_{i=0}^d; A;  \{E_{d-i}\}_{i=0}^d)$ \\ 
$\Phi^{\downarrow \Downarrow *}$ &
       $(A^*; \{E^*_{d-i}\}_{i=0}^d; A;  \{E_{d-i}\}_{i=0}^d)$
\end{tabular}
\end{center}

\begin{definition}
\rm
Let $\Phi$ denote a tridiagonal system.
For $g \in D_4$ and for an object $f$ associated with $\Phi$, we let
$f^g$ denote the corresponding object associated with $\Phi^{g^{-1}}$
(we have been using this convention all along; an example is
$\theta^*_i(\Phi)=
\theta_i(\Phi^*)$).
We say $f$ is {\it $D_4$-invariant} whenever $f^g=f$
for all $g \in D_4$. 
\end{definition}

\noindent Let $\Phi$ denote a tridiagonal system over $\K$.
In this paper we associate with $\Phi$ a certain polynomial
$P$ called the Drinfel'd polynomial,
and show 
that $P$ is $D_4$-invariant.

\medskip
\noindent For later use we remark that 
the elements $*, \Downarrow$ together generate $D_4$.

\section{The eigenvalues and dual eigenvalues}

\noindent In order to develop our theory of the Drinfel'd 
polynomial, we will need some detailed supporting
results concerning three sequences of scalars:
the eigenvalue sequence, the dual eigenvalue sequence,
and the split sequence. The supporting results
are contained in this section and the next two.

\begin{definition}        \label{def}
\rm
Let $\Phi=(A;\{E_i\}_{i=0}^d;A^*$; $\{E^*_i\}_{i=0}^d)$ 
denote a tridiagonal system on $V$.
For $0 \leq i \leq d$ let $\theta_i$ (resp. $\theta^*_i$)
denote the eigenvalue of $A$ (resp. $A^*$)
associated with the eigenspace $E_iV$ (resp. $E^*_iV$).
We call $\{\theta_i\}_{i=0}^d$ (resp. $\{\theta^*_i\}_{i=0}^d$)
the {\em eigenvalue sequence}
(resp. {\em dual eigenvalue sequence}) of $\Phi$.
We observe that $\{\theta_i\}_{i=0}^d$ (resp. $\{\theta^*_i\}_{i=0}^d$) are
mutually distinct and contained in $\K$.
\end{definition}
\medskip

\begin{lemma} {\rm \cite[Theorem 11.1]{TD00}}  \label{lem:indep} \samepage
With reference to Definition \ref{def}, the expressions
\begin{equation}         \label{eq:beta}
 \frac{\theta_{i-2}-\theta_{i+1}}{\theta_{i-1}-\theta_i},
\qquad\qquad
 \frac{\theta^*_{i-2}-\theta^*_{i+1}}{\theta^*_{i-1}-\theta^*_i}
\end{equation}
are equal and independent of $i$ for $2 \leq i \leq d-1$.
\end{lemma}

\begin{definition}
\label{def:base}
\rm
Let $A,A^*$ denote a tridiagonal pair over $\K$.
We associate with $A,A^*$ 
a scalar $\beta \in \K$
as follows. 
If the diameter  $d\geq 3$
let $\beta+1$ denote the common value of
(\ref{eq:beta}), where 
 $\{\theta_i\}_{i=0}^d$ (resp. 
 $\{\theta^*_i\}_{i=0}^d$) is a standard ordering
of the eigenvalues of $A$ (resp. $A^*$).
If $d \leq 2$ let $\beta$ 
denote any nonzero scalar in $\K$. 
We call $\beta$ the {\it base} of $A,A^*$.
By construction, for $d\geq 3$ 
the tridiagonal pairs
$A,A^*$ and $A^*, A$ have the same base.
For $d\leq 2$, 
we always choose the bases such that 
$A,A^*$ and $A^*, A$ have the same base.
\end{definition}

\begin{definition}
\label{def:type}
\rm
Let $A,A^*$ denote a tridiagonal pair over $\K$,
with diameter $d$ and base $\beta$.
We assign to $A,A^*$ a {\it type} as follows:
\begin{center}
\begin{tabular}{c|c}
type  &  description\\
\hline
${\rm I}$ & 
       $\beta \not=2$, $\beta \not=-2$ \\ 
${\rm II}$ &
       $\beta=2$, ${\rm Char}(\K)\not=2$ \\ 
${\rm III}^+$ &
       $\beta=-2$, ${\rm Char}(\K)\not=2$, $d$ even \\ 
${\rm III}^-$ &
       $\beta=-2$, ${\rm Char}(\K)\not=2$, $d$ odd \\ 
${\rm IV}$ &
       $\beta=0$, ${\rm Char}(\K)=2$ 
\end{tabular}
\end{center}
By 
{\rm  \cite[Theorem~11.2]{TD00}},
if $A,A^*$ has type IV then $d=3$.
We say $A,A^*$ has {\it type} III whenever $A,A^*$ has type
${\rm III}^+$ or
${\rm III}^-$.
\end{definition}

\begin{definition}   
     \label{defe}
\rm
Let $\Phi$
denote a tridiagonal system.
By the {\it base} (resp. {\it type}) of $\Phi$ we mean the base
(resp. type) of the associated tridiagonal pair.
By construction, the base (resp. type) of $\Phi$ is
$D_4$-invariant.
\end{definition}
\medskip

\noindent
Let $\overline \K$ denote the
algebraic closure of $\K$.

\begin{lemma}
{\rm  \cite[Theorem~11.2]{TD00}}
\label{lem:eig1}
With reference to Definition \ref{def},
 assume
$\Phi$ has type I,
and fix a nonzero $q \in {\overline \K}$ 
such that $q^2+q^{-2}=\beta$. Then there exists 
a sequence of scalars 
$a,b,c,a^*,b^*,c^*$ taken from ${\overline \K}$ such that
\begin{eqnarray*}
\theta_i &=& a + b q^{2i-d}+c q^{d-2i}, 
\\
\theta^*_i &=& a^* +  b^* q^{2i-d}+c^* q^{d-2i} 
\end{eqnarray*}
for $0 \leq i \leq d$.
The sequence is uniquely determined by $q$ provided $d\geq 2$.
\end{lemma}

\begin{lemma}
\label{lem:l2I}
With reference to Definition 
 \ref{def}
and Lemma \ref{lem:eig1}, for $0 \leq i,j\leq d$ we have
\begin{eqnarray*}
\theta_i - \theta_j &=& (q^{i-j}-q^{j-i})(bq^{i+j-d}-cq^{d-i-j}), \\
\theta^*_i - \theta^*_j &=& (q^{i-j}-q^{j-i})(b^*q^{i+j-d}-c^*q^{d-i-j}).
\end{eqnarray*}
\end{lemma}

\begin{note}
\label{lem:l3I}
\rm
With reference to Definition 
 \ref{def}
and Lemma 
\ref{lem:eig1}, for $1 \leq i \leq d$ we have $q^{2i}\not=1$;
otherwise $\theta_i=\theta_0$ by Lemma
\ref{lem:l2I}.
\end{note}

\begin{lemma}
\label{lem:oioiI}
With reference to Definition 
 \ref{def}
and Lemma \ref{lem:eig1}, for $1 \leq i \leq d$
we have
\begin{eqnarray*}
\frac{(\theta_0-\theta_i)
(\theta^*_0-\theta^*_i)}{
(q^i-q^{-i})^2} = 
bb^*q^{2i-2d}+cc^*q^{2d-2i}-bc^*-cb^*.
\end{eqnarray*}
\end{lemma}
\noindent {\it Proof:} Use
Lemma
\ref{lem:l2I}.
\hfill $\Box$ \\

\begin{lemma}
{\rm  \cite[Theorem~11.2]{TD00}}
\label{lem:eig2}
With reference to Definition \ref{def},
assume
$\Phi$ has type II.
Then there exists 
a sequence of scalars 
$a,b,c,a^*,b^*,c^*$ taken from $\K$ such that
\begin{eqnarray*}
\theta_i &=& a+ b(i-d/2)+ci(d-i), 
\\
\theta^*_i &=& a^*+ b^*(i-d/2)+c^*i(d-i)
\end{eqnarray*}
for $0 \leq i \leq d$.
The sequence is unique provided $d\geq 2$.
\end{lemma}

\begin{lemma}
\label{lem:l2}
With reference to Definition \ref{def}
and 
 Lemma \ref{lem:eig2}, for $0 \leq i,j\leq d$ we have
\begin{eqnarray*}
\theta_i - \theta_j &=& (i-j)(b+c(d-i-j)), \\
\theta^*_i - \theta^*_j &=& (i-j)(b^*+c^*(d-i-j)).
\end{eqnarray*}
\end{lemma}

\begin{note}
\label{lem:l3II}
\rm
With reference to Definition \ref{def}
and Lemma \ref{lem:eig2},
for all primes $p\leq d$
we have $\mbox{\rm Char}(\K)\not=p$; otherwise
$\theta_p=\theta_0$ by Lemma  
\ref{lem:l2}.
Consequently 
$\mbox{\rm Char}(\K)$ is 0 or an odd prime greater than $d$.
\end{note}

\begin{note}
\label{lem:l3IIa}
\rm
With reference to Definition \ref{def}
and Lemma \ref{lem:eig2}, assume $d\geq 1$. Then $b \not=0$; otherwise
$\theta_0=\theta_d$. Similary $b^*\not=0$.
\end{note}

\begin{lemma}
\label{lem:oioiII}
With reference to Definition \ref{def}
and Lemma 
 \ref{lem:eig2}, for $1 \leq i \leq d$
we have
\begin{eqnarray*}
(\theta_0-\theta_i)
(\theta^*_0-\theta^*_i)i^{-2} = 
b b^* +(b c^* +c b^*)(d-i) + c c^* (d-i)^2.
\end{eqnarray*}
\end{lemma}
\noindent {\it Proof:}
Use Lemma \ref{lem:l2}.
\hfill $\Box$ \\

\begin{lemma}
{\rm  \cite[Theorem~11.2]{TD00}}
\label{lem:eig3}
With reference to Definition \ref{def},
assume
$\Phi$ has type III.
Then there exists 
a sequence of scalars 
$a,b,c,a^*,b^*,c^*$ taken from ${\K}$ such that
\begin{eqnarray*}
\theta_i &=&     
\begin{cases}  
a + b +c(i-d/2)        
 & \text{\rm if $i$ is even},  \\
a - b -c(i-d/2)        
 & \text{\rm if $i$ is odd},
    \end{cases} 
\\
\theta^*_i &=&     
\begin{cases}  
a^* + b^* +c^*(i-d/2)        
 & \text{\rm if $i$ is even},  \\
a^* - b^* -c^*(i-d/2)        
 & \text{\rm if $i$ is odd}
    \end{cases}
\end{eqnarray*}
for $0 \leq i \leq d$.
The sequence is unique provided $d\geq 2$.
\end{lemma}

\begin{lemma}
\label{lem:l3}
With reference to Definition \ref{def}
and Lemma 
 \ref{lem:eig3}, for $0 \leq i,j\leq d$ we have
\begin{eqnarray*}
\theta_i - \theta_j  &=& 
\begin{cases}  
(-1)^i c(i-j)        
 & \text{\rm if $i-j$ is even},  \\
(-1)^i(2b+c(i+j-d))        
& \text{\rm if $i-j$ is odd},
    \end{cases}
\\
\theta^*_i - \theta^*_j  &=& 
\begin{cases}  
(-1)^i c^*(i-j)        
 & \text{\rm if $i-j$ is even},  \\
(-1)^i(2b^*+c^*(i+j-d))        
& \text{\rm if $i-j$ is odd}.
    \end{cases}
\end{eqnarray*}
\end{lemma}

\begin{note}
\label{lem:l3III}
\rm
With reference to Definition \ref{def}
and Lemma \ref{lem:eig3}, 
for all primes $p \leq d/2 $
  we have 
$\mbox{\rm Char}(\K)\not=p$; otherwise
$\theta_{2p}=\theta_0$
by Lemma 
\ref{lem:l3}.
 Consequently 
$\mbox{\rm Char}(\K)$ is 0 or an odd prime greater than $d/2$.
\end{note}

\begin{note}
\label{lem:l3IIIa}
\rm
With reference to Definition \ref{def}
and Lemma \ref{lem:eig3},
assume $d\geq 2$. Then 
$c \not=0$; otherwise
$\theta_0=\theta_2$. Similarly $c^*\not=0$.
For $d\leq 1$ we always choose $c,c^*$ to be nonzero.
\end{note}

\begin{note}
\label{lem:l3IIIab}
\rm
With reference to Definition \ref{def}
and Lemma \ref{lem:eig3},
assume $d$ is odd. Then $b \not=0$;
otherwise $\theta_0=\theta_d$.
Similarly $b^*\not=0$.
\end{note}

\begin{lemma}
\label{lem:oioiIIIp}
With reference to Definition \ref{def}
and Lemma \ref{lem:eig3}, for $0 \leq i \leq d$
we have
\begin{eqnarray*}
(\theta_0-\theta_i)
(\theta^*_0-\theta^*_i) = 
\begin{cases}  
cc^*i^2        
 & \text{\rm if $i$ is even},  \\
4b b^* +2(b c^* + c b^*)(i-d)+ c c^* (i-d)^2
& \text{\rm if $i$ is odd}.
    \end{cases}
\end{eqnarray*}
\end{lemma}
\noindent {\it Proof:}
Use Lemma
\ref{lem:l3}.
\hfill $\Box$ \\

\begin{lemma}
{\rm  \cite[Theorem~11.2]{TD00}}
\label{lem:eig4}
With reference to Definition \ref{def},
 assume
$\Phi$ has type IV.
Then there exists 
a unique sequence of scalars 
$a,b,c,a^*,b^*,c^*$ taken from ${\K}$ such that
\begin{eqnarray*}
&&
\theta_0=a,
\qquad 
\theta_1=b+c,
\qquad 
\theta_2=a+c,
\qquad 
\theta_3=b,
\\
&&
\theta^*_0=a^*,
\qquad 
\theta^*_1=b^*+c^*,
\qquad 
\theta^*_2=a^*+c^*,
\qquad 
\theta^*_3=b^*.
\end{eqnarray*}
\end{lemma}

\begin{note}
\label{lem:l3IV}
\rm
With reference to Definition \ref{def}
and Lemma \ref{lem:eig4},
each of $a+b, a+b+c, c$ is nonzero since
$\lbrace\theta_i\rbrace_{i=0}^3$ are distinct.
Similarly
each of $a^*+b^*, a^*+b^*+c^*, c^*$ is nonzero.
\end{note}

\begin{lemma}
\label{lem:oioiIV}
With reference to Definition \ref{def}
and Lemma \ref{lem:eig4},
\begin{eqnarray*}
(\theta_0-\theta_1)
(\theta^*_0-\theta^*_1) &=& (a+b+c)(a^*+b^*+c^*),
\\
(\theta_0-\theta_2)
(\theta^*_0-\theta^*_2) &=& cc^*,
\\
(\theta_0-\theta_3)
(\theta^*_0-\theta^*_3) &=& (a+b)(a^*+b^*).
\end{eqnarray*}
\end{lemma}
\noindent {\it Proof:}
Use Lemma \ref{lem:eig4}.
\hfill $\Box$ \\

\medskip
\noindent
When discussing tridiagonal systems we will use the following
notational convention.
Let $\lambda$ denote an indeterminate and let $\K[\lambda]$
denote the $\K$-algebra consisting of the polynomials
in $\lambda$ that have all coefficients in $\K$.
With reference to Definition \ref{def},
for $0 \leq i \leq d$ we define the following polynomials in
$\K[\lambda]$:
\begin{eqnarray*}
 \tau_i &=& 
  (\lambda-\theta_0)(\lambda-\theta_1)\cdots(\lambda -\theta_{i-1}), \\
 \eta_i &=&
  (\lambda-\theta_d)(\lambda-\theta_{d-1})\cdots(\lambda-\theta_{d-i+1}),  \\
 \tau^*_i &=&
  (\lambda-\theta^*_0)(\lambda-\theta^*_1)\cdots(\lambda-\theta^*_{i-1}), \\
 \eta^*_i &=&
  (\lambda-\theta^*_d)(\lambda-\theta^*_{d-1})\cdots(\lambda-\theta^*_{d-i+1}).
\end{eqnarray*}
Note that each of $\tau_i$, $\eta_i$, $\tau^*_i$, $\eta^*_i$ is
monic with degree $i$.

\begin{lemma}
\label{lem:tauzero}
With reference to Definition \ref{def},
the following 
{\rm (i)}--{\rm (iv)} 
hold for $0 \leq i,j\leq d$.
\begin{itemize}
\item[\rm (i)]
$\tau_i(\theta_j)=0$ if and only if $j<i$;
\item[\rm (ii)]
$\eta_i(\theta_j)=0$ if and only if $j>d-i$;
\item[\rm (iii)]
$\tau^*_i(\theta^*_j)=0$ if and only if $j<i$;
\item[\rm (iv)]
$\eta^*_i(\theta^*_j)=0$ if and only if $j>d-i$.
\end{itemize}
\end{lemma}

\section{Some scalars}

\noindent 
Adopt the notation of Definition \ref{def}.  
For nonnegative integers $r,s,t$ such that $r+s+t\leq d$,
in
   \cite[Definition 13.1]{LS24}
we defined some scalars 
 $\lbrack r,s,t\rbrack_q \in \K$.
   By \cite[Definition 13.1]{LS24}
 these scalars
are rational functions of the base $\beta$,
and in this paper we are going to drop the
subscript $q$ altogether and just write 
 $\lbrack r,s,t\rbrack$.
These scalars are described in the  next definition.
We will use the following notation.
For all $a,q \in {\overline \K}$ define
\begin{eqnarray}
(a;q)_n = (1-a)(1-aq)\cdots (1-aq^{n-1}) \qquad  \qquad n=0,1,2\ldots
\label{eq:aqn}
\end{eqnarray}
We interpret 
$(a;q)_0=1$.

\begin{definition}
{\rm \cite[Lemma 13.2]{LS24}}
\label{lem:brackform}
\rm
With reference to Definition \ref{def},
let $r,s,t$ denote nonnegative integers such that $r+s+t\leq d$.
We define $\lbrack r,s,t \rbrack$ as follows.
\begin{itemize}
\item[\rm (i)]
Assume $\Phi$ is type I. Then
\begin{eqnarray}
 \lbrack r,s,t\rbrack= 
\frac{
(q^2;q^2)_{r+s}
(q^2;q^2)_{r+t}
(q^2;q^2)_{s+t}
}
{
(q^2;q^2)_{r}
(q^2;q^2)_{s}
(q^2;q^2)_{t}
(q^2;q^2)_{r+s+t}
}.
\label{eq:rst}
\end{eqnarray}
Here $q^2+q^{-2}=\beta$ where $\beta$ is the base of $\Phi$. 
\item[\rm (ii)] 
Assume $\Phi$ is type II.  Then
\begin{equation}
\lbrack r,s,t \rbrack= {{(r+s)!\,(r+t)!\,(s+t)!}\over{r!\,s!\,t!\,(r+s+t)!}}.
\label{eq:nextlevelone}
\end{equation}
\item[\rm (iii)] 
Assume $\Phi$ is type III. 
If each of $r,s,t$ is odd, then
$\lbrack r,s,t\rbrack =0$. If at least one of $r,s,t$ is even,
then
\begin{equation}
\lbrack r,s,t \rbrack = 
\frac{\lfloor \frac{r+s}{2}\rfloor ! 
\lfloor \frac{r+t}{2}\rfloor !
\lfloor \frac{
s+t}{2}\rfloor !}{\lfloor \frac{r}{2}\rfloor ! 
\lfloor \frac{s}{2}\rfloor ! 
\lfloor \frac{t}{2}\rfloor !  \lfloor \frac{r+s+t}{2} \rfloor !} .
\label{eq:nextlevelminone}
\end{equation}
The expression $\lfloor x \rfloor $ denotes the greatest integer less
than or equal to $x$.
\item[\rm (iv)] 
Assume $\Phi$ is type IV. 
Recall in this case  $d=3$.
If each of $r,s,t$ equals 1, then $\lbrack r,s,t\rbrack =0$.
If at least one of $r,s,t$ equals 0 then $\lbrack r,s,t \rbrack=1$.
\end{itemize}
\end{definition}

\noindent We have a comment.
\begin{lemma}
\label{lem:sym} 
With reference to Definitions \ref{def}
and 
\ref{lem:brackform}
the following
{\rm (i)}, {\rm (ii)} hold.
\begin{itemize}
\item[\rm (i)]
 The expression $\lbrack r,s,t\rbrack$
is symmetric in $r,s,t$. 
\item[\rm (ii)] If at least one of $r,s,t$  is zero then
 $\lbrack r,s,t\rbrack=1$.
\item[\rm (iii)]
The expression 
  $\lbrack r,s,t\rbrack$ is $D_4$-invariant.
\end{itemize}
\end{lemma}

\begin{lemma}
\label{lem:double}
With reference to Definition \ref{def},
let $r,s,t,u$ denote nonnegative integers whose sum is at most $d$. Then
\begin{eqnarray*}
\lbrack r,s,t+u\rbrack
\;\lbrack t,u,r+s\rbrack
=
\lbrack s,u,r+t\rbrack
\;\lbrack r,t,s+u\rbrack.
\end{eqnarray*}
\end{lemma}
\noindent {\it Proof:}
For each type I--IV this is a routine cancellation
using the formulae in 
Definition \ref{lem:brackform}.
\hfill $\Box$ \\

\noindent The following lemma shows one significance of
the scalars $\lbrack r,s,t\rbrack$.

\begin{lemma} 
\label{lem:etatau}
\rm
\cite[Lemma 9.1]{nomsharp}
With reference to Definition \ref{def},
\begin{eqnarray*}
\eta_i = \sum_{h=0}^i \lbrack h,i-h,d-i\rbrack \eta_{i-h}
(\theta_0)
 \tau_h
\qquad \qquad (0 \leq i \leq d).
\end{eqnarray*}
\end{lemma}

\section{The split decomposition}
\indent
In this section we recall the split decomposition associated with
a tridiagonal system \cite[Section 4]{TD00}.
With reference to Definition \ref{def},
for $0 \leq i \leq d$ define
\begin{equation}          \label{eq:defUi}
   U_i = (E^*_0V+E^*_1V+\cdots+E^*_iV) \cap(E_iV+E_{i+1}V+\cdots+E_dV).
\end{equation}
By \cite[Theorem 4.6]{TD00}
\[
   V = U_0+U_1+\cdots+U_d  \qquad\qquad (\mbox{\rm direct sum}),
\]
and for $0 \leq i \leq d$ both
\begin{eqnarray}
  U_0+U_1+\cdots+U_i &=& E^*_0V+E^*_1V+\cdots+E^*_iV,    \label{eq:sumU0Ui} \\
  U_i+U_{i+1}+\cdots+U_d &=& E_iV+E_{i+1}V+\cdots+E_dV.  \label{eq:sumUiUd}
\end{eqnarray}
By \cite[Corollary 5.7]{TD00}  $U_i$ has dimension $\rho_i$
for $0 \leq i \leq d$, where $\{\rho_i\}_{i=0}^d$ is the shape of $\Phi$.
By \cite[Theorem 4.6]{TD00} both
\begin{eqnarray}
  (A-\theta_i I)U_i & \subseteq & U_{i+1},    \label{eq:up}  \\
  (A^*-\theta^*_i I)U_i & \subseteq & U_{i-1}  \label{eq:down}
\end{eqnarray}
for $0 \leq i \leq d$, where $U_{-1}=0$ and $U_{d+1}=0$.
The sequence $\{U_i\}_{i=0}^d$ is called the {\em $\Phi$-split decomposition}
of $V$ \cite[Section 4]{TD00}.
Now assume that $\Phi$ is sharp, so 
that $U_0$ has dimension 1.
By 
  (\ref{eq:up}),
  (\ref{eq:down}),
for $0 \leq i \leq d$ the space $U_0$ is invariant under
\begin{equation}        \label{eq:updown}
 (A^*-\theta^*_1I)(A^*-\theta^*_2I)\cdots(A^*-\theta^*_iI)
 (A-\theta_{i-1}I)\cdots(A-\theta_1I)(A-\theta_0I);
\end{equation}
let $\zeta_i$ denote the corresponding eigenvalue.
Note that $\zeta_0=1$.
We call the sequence $\{\zeta_i\}_{i=0}^d$ the {\em split sequence}
of $\Phi$.

\begin{note}
\label{nt:spltdiff}
\rm
In the literature on Leonard systems there
are two sequences of scalars
called the first split sequence and the
second split sequence \cite[Section~3]{LS99}. 
These sequences are related to
the above split sequence as follows. 
With reference to
Definition
 \ref{def}, assume $\Phi$ is a Leonard system and
let $\lbrace \varphi_i\rbrace_{i=1}^d$ 
(resp. 
 $\lbrace \phi_i\rbrace_{i=1}^d$) 
denote the
first split sequence (resp. second split sequence)
for $\Phi$ in the sense of \cite[Section~3]{LS99}.
Then the sequence
 $\lbrace \varphi_1\varphi_2\cdots \varphi_i\rbrace_{i=0}^d$
(resp. 
 $\lbrace \phi_1\phi_2\cdots \phi_i\rbrace_{i=0}^d$)
is the split sequence of $\Phi$ (resp.
$\Phi^\Downarrow$).
\end{note}

\noindent The $D_4$ action affects the
split sequence as follows.

\begin{lemma}
\label{lem:dualsp}
{\rm \cite[Theorem 7.3]{nomsharp}}
With reference to Definition \ref{def} assume
$\Phi$ is sharp. Then
\begin{eqnarray*}
 \zeta^*_i=\zeta_i 
\qquad \qquad 
0 \leq i \leq d.
\end{eqnarray*}
\end{lemma}

\begin{lemma}
\label{lem:zetadown}
{\rm \cite[Theorem 9.3]{nomsharp}}
With reference to Definition \ref{def} assume
$\Phi$ is sharp. Then 
\begin{eqnarray*}
\zeta^{\Downarrow}_i = 
\sum_{h=0}^i \lbrack h,i-h,d-i\rbrack
 \frac{\eta^*_{d-h}(\theta^*_0)  
\eta_{i-h}(\theta_0) 
}
{
 \eta^*_{d-i}(\theta^*_0)
}
\; \zeta_h
\qquad \qquad  0 \leq i \leq d.
\end{eqnarray*}
\end{lemma}

\begin{definition}
{\rm \cite[Definition~6.2]{nomsharp}}
\rm
Let $\Phi$
denote a sharp tridiagonal system.
By the {\it parameter array} of $\Phi$ we mean the
sequence
 $(\{\theta_i\}_{i=0}^d; \{\theta^*_i\}_{i=0}^d; \{\zeta_i\}_{i=0}^d)$
where
 $\{\theta_i\}_{i=0}^d$
(resp.
$\{\theta^*_i\}_{i=0}^d$
)
is the eigenvalue sequence
(resp. dual eigenvalue sequence)
of $\Phi$ and
$\{\zeta_i\}_{i=0}^d$ is the split sequence of $\Phi$.
\end{definition}

\begin{proposition}
{\rm \cite[Theorem~1.6]{nomstructure}}
Two sharp tridiagonal systems over $\K$ are isomorphic
if and only if they have the same parameter array.
\end{proposition}

\section{The Drinfel'd polynomial}

\noindent
Let $\Phi$ denote a sharp tridiagonal system.
In this section we introduce the
Drinfel'd polynomial of $\Phi$, for all types
except 
 ${\rm III}^{+}$.
In Section 8 we will 
treat type
 ${\rm III}^{+}$.

\begin{definition}
\label{def:aap}
\rm
Let $\Phi$ denote a sharp tridiagonal system over 
$\K$ with eigenvalue sequence
 $\lbrace \theta_i\rbrace_{i=0}^d$
and dual eigenvalue sequence  
 $\lbrace \theta^*_i\rbrace_{i=0}^d$.
Assume $\Phi$ is not type
 ${\rm III}^{+}$.
For $1 \leq i \leq d$ 
we define a polynomial $p_i \in \F\lbrack \lambda \rbrack$
by 
\begin{eqnarray}
\label{eq:drinprelim}
p_i = (\theta_0 \theta^*_d+\theta_d \theta^*_0-\lambda)\alpha_i
+
(\theta_0-\theta_i) 
(\theta^*_0-\theta^*_i),
\end{eqnarray}
where 
\begin{center}
\begin{tabular}{c|c c c c c c}
Case  &  I & II & ${\rm III}^{-}$($i$ even) &
 ${\rm III}^{-}$($i$ odd) &
 IV ($i$ even) & IV ($i$ odd) \\
\hline
$\alpha_i$ & 
$\frac{(q^i-q^{-i})^2
}{(q^d-q^{-d})^2}$
&      
$\frac{i^2
}{d^2}$
&
0
& 
$1$
& 
0
&
1
\end{tabular}
\end{center}
For type {\rm I},
 $q^2+q^{-2}=\beta$ where
 $\beta$ is the base of $\Phi$ from Definition
     \ref{defe}.
\end{definition}

\begin{note}
\rm Referring to the table in Definition
\ref{def:aap}, for $1 \leq i \leq d$
the expression $\alpha_i$ is a rational function
of the base $\beta$. For example if $d=3$,
\begin{eqnarray*}
\alpha_1 = (\beta+1)^{-2},
\qquad \qquad 
\alpha_2 = (\beta+2)(\beta+1)^{-2},
\qquad \qquad 
\alpha_3 =1.
\end{eqnarray*}
\end{note}

\begin{lemma}
\label{lem:drin0}
With reference to
Definition
\ref{def:aap}
the following 
 {\rm (i)}, {\rm (ii)}
hold for 
$1 \leq i \leq d$.
\begin{itemize}
\item[\rm (i)]
 $\alpha_i$
is $D_4$-invariant. 
\item[\rm (ii)]
$p^*_i=p_i$.
\end{itemize}
\end{lemma}
\noindent {\it Proof:}
(i) Follows from
the table
in Definition
\ref{def:aap}, and since $\beta$ is $D_4$-invariant.
\\
\noindent (ii) Immediate from
(\ref{eq:drinprelim}) and (i) above.
\hfill $\Box$ \\

\begin{lemma}
\label{def:aap2}
With reference to Definition \ref{def:aap},
for $1 \leq i \leq d$ both
\begin{eqnarray*}
p_i(
\theta_0 \theta^*_d+\theta_d \theta^*_0)
 &=& (\theta_0-\theta_i)(\theta^*_0-\theta^*_i),
\\
p_i(
\theta_0 \theta^*_0+\theta_d \theta^*_d)
 &=& (\theta_0-\theta_i)(\theta^*_0-\theta^*_i)-\alpha_i
 (\theta_0-\theta_d)(\theta^*_0-\theta^*_d).
\end{eqnarray*}
\end{lemma}

\begin{lemma}
\label{lem:aapd}
With reference to 
Definition \ref{def:aap},
\begin{eqnarray}
\label{eq:piddown}
p^\Downarrow_i = (\theta_0 \theta^*_0+\theta_d \theta^*_d-\lambda)\alpha_i
+
(\theta_d-\theta_{d-i}) 
(\theta^*_0-\theta^*_i)
\qquad \qquad (1 \leq i \leq d).
\end{eqnarray}
\end{lemma}
\noindent {\it Proof:}
Apply $\Downarrow$ to 
(\ref{eq:drinprelim})
and evaluate the result using
Lemma
\ref{lem:drin0}(i) 
and 
$\theta^{\Downarrow}_j = \theta_{d-j}$
 $(0 \leq j \leq d)$.
\hfill $\Box$ \\

\begin{lemma}
\label{def:aapd2}
With reference to Definition \ref{def:aap},
for $1 \leq i \leq d$ both
\begin{eqnarray*}
p^\Downarrow_i(
\theta_0 \theta^*_0+\theta_d \theta^*_d)
 &=& (\theta_d-\theta_{d-i})(\theta^*_0-\theta^*_i),
\\
p^\Downarrow_i(
\theta_0 \theta^*_d+\theta_d \theta^*_0)
 &=& (\theta_d-\theta_{d-i})(\theta^*_0-\theta^*_i)+\alpha_i
 (\theta_0-\theta_d)(\theta^*_0-\theta^*_d).
\end{eqnarray*}
\end{lemma}

\begin{lemma}
\label{lem:dnorm}
With reference to
Definition \ref{def:aap}, for $d\geq 1$ both
\begin{eqnarray*}
p_d &=& \theta_0 \theta^*_0 + \theta_d \theta^*_d - \lambda,
\\
p^{\Downarrow}_d &=& \theta_0 \theta^*_d + \theta_d \theta^*_0 - \lambda.
\end{eqnarray*}
\end{lemma}
\noindent {\it Proof:}
Set $i=d$ in
(\ref{eq:drinprelim}), 
(\ref{eq:piddown})
and observe that $\alpha_d=1$.
\hfill $\Box$ \\

\begin{definition}
\label{def:pq}
\rm
With reference to Definition
\ref{def:aap},
we define a polynomial $P \in \F\lbrack \lambda \rbrack$
by 
\begin{eqnarray}
P = \sum_{i=0}^d \zeta_i p_{i+1}p_{i+2}\cdots p_d,
\label{eq:drinp}
\end{eqnarray}
where $\lbrace \zeta_i\rbrace_{i=0}^d$ is the split
sequence of $\Phi$ from Section 6.
We call $P$ the {\it Drinfel'd polynomial} of $\Phi$.
\end{definition}

\begin{example}
\rm
With reference to Definition \ref{def:aap} and
Definition 
\ref{def:pq}, if $d=0$ then $P=1$, and if $d=1$ then
$P=\zeta_1 + \theta_0 \theta^*_0 + \theta_1 \theta^*_1 - \lambda$.
\end{example}

\noindent 
With reference to Definition \ref{def:aap} and
Definition 
\ref{def:pq}, we are going to show that
$P$ is $D_4$-invariant.
 We will prove this after a few lemmas.

\begin{lemma}
\label{lem:dualsame}
With reference to Definitions \ref{def:aap} and
\ref{def:pq}, we have
$P^*=P$.
\end{lemma}
\noindent {\it Proof:}
Apply $*$ to 
(\ref{eq:drinp}) 
and evaluate the result using
Lemma \ref{lem:dualsp} and
Lemma \ref{lem:drin0}(ii).
\hfill $\Box$ \\

\begin{lemma}
\label{lem:cijpre}
With reference to Definition
\ref{def:aap},
for $1 \leq i\leq j \leq d$ we have
\begin{eqnarray}
\label{eq:alal}
\lbrack i,j-i,d-j\rbrack \alpha_i = 
\lbrack i-1,j-i,d-j+1\rbrack \alpha_j a_{i,j}
\end{eqnarray}
where $a_{i,j}$ is given below.
{\rm
\begin{center}
\begin{tabular}{c|c c c}
type of $\Phi$  &  I & II & 
 IV  \\
\hline
$a_{i,j}$ & 
$\frac{(q^i-q^{-i})(
q^{d-j+1}-q^{j-d-1})}
{
(q^j-q^{-j})
(q^{d-i+1}-q^{i-d-1})
}$
&      
$\frac{i(d-j+1)
}{j(d-i+1)}$
&
$1$
\end{tabular}
\end{center}
For type ${\rm III}^{-}$ the scalar
$a_{i,j}$ depends on the parity
of $i,j$ as follows:
\begin{center}
\begin{tabular}{c|c c c c}
Parity of $i,j$  & $i$ even, $j$ even & $i$ odd, $j$ even & 
$i$ even, $j$ odd & $i$ odd, $j$ odd  
  \\
\hline
$a_{i,j}$ & 
$\frac{i(d-j+1)
}
{j(d-i+1)
}$
&   
$
\frac{j-d-1}
{j}$
&   
$\frac{i}
{i-d-1}
$
&   
$1$
\end{tabular}
\end{center}
}
\end{lemma}
\noindent {\it Proof:}
For each of the above seven entries
one routinely verifies
(\ref{eq:alal}) using
Definition \ref{lem:brackform}
and the table in
Definition
\ref{def:aap}.
\hfill $\Box$ \\

\begin{lemma}
\label{lem:cij}
With reference to Definition
\ref{def:aap},
for $1 \leq i\leq j\leq d$ we have
\begin{eqnarray}
\label{eq:ppc}
\lbrack i,j-i,d-j\rbrack
p^\Downarrow_i &=& 
\lbrack i-1,j-i,d-j+1\rbrack a_{i,j} p_j
\;+\; c_{i,j}
\end{eqnarray}
where
$a_{i,j}$ is from
Lemma
\ref{lem:cijpre}
and 
\begin{eqnarray*}
&&c_{i,j}\;=\; 
\lbrack i,j-i,d-j\rbrack (\theta_d -\theta_{d-i})(\theta^*_0-\theta^*_i)
\\
&&
\qquad \qquad \qquad - \quad
\lbrack i-1,j-i,d-j+1\rbrack a_{i,j} p_j(\theta_0 \theta^*_0
+\theta_d \theta^*_d).
\label{eq:cijform}
\end{eqnarray*}
\end{lemma}
\noindent {\it Proof:}
Each side of 
(\ref{eq:ppc}) is a polynomial in 
$\K\lbrack \lambda  \rbrack$
with degree at most 1.
On each side of
(\ref{eq:ppc}) 
the
$\lambda$-coefficients agree
by
(\ref{eq:drinprelim}), 
(\ref{eq:piddown}) and 
 Lemma
\ref{lem:cijpre}.
To see that the constant terms also agree,
evaluate
each side of
(\ref{eq:ppc}) 
at $\lambda=
\theta_0\theta^*_0
+
\theta_d\theta^*_d
$
and use the first equation of
Lemma
\ref{def:aapd2}.
\hfill $\Box$ \\

\begin{lemma}
\label{prop:aapdvsp}
With reference to
Definition
\ref{def:aap},
for $0 \leq i \leq d$ we have
\begin{eqnarray*}
p^\Downarrow_{i+1}
p^\Downarrow_{i+2}
\cdots
p^\Downarrow_{d}
&=& 
\sum_{h=i}^d 
\lbrack i,h-i,d-h\rbrack
 \frac{\eta^*_{d-i}(\theta^*_0)  
\tau_{h-i}(\theta_d) 
}
{
 \eta^*_{d-h}(\theta^*_0)  
}
p_{h+1}
p_{h+2}
\cdots
p_{d}.
\end{eqnarray*}
\end{lemma}
\noindent {\it Proof:}
The proof is by induction on $i=d,d-1,\ldots, 0$.
Let $i$ be given. The assertion for $i=d$ holds trivially,
so assume $i<d$. By induction,
\begin{eqnarray*}
p^\Downarrow_{i+2}
\cdots
p^\Downarrow_{d}
&=&
\sum_{h=i+1}^d 
\lbrack i+1,h-i-1,d-h\rbrack
 \frac{\eta^*_{d-i-1}(\theta^*_0)  
\tau_{h-i-1}(\theta_d) 
}
{
 \eta^*_{d-h}(\theta^*_0)  
}
p_{h+1}
\cdots
p_{d}.
\end{eqnarray*}
In this equation we multiply both sides
by $p^\Downarrow_{i+1}$
and use 
Lemma
\ref{lem:cij}
to get
\begin{eqnarray*}
&& p^\Downarrow_{i+1}
p^\Downarrow_{i+2}
\cdots
p^\Downarrow_{d}
\quad  = \quad
\sum_{h=i+1}^d 
\Bigl(
\lbrack i,h-i-1,d-h+1\rbrack a_{i+1,h}p_h
+c_{i+1,h}
\Bigr)
\\
&&
 \qquad \qquad  
 \qquad \qquad  \qquad  \qquad 
\times \quad 
 \frac{\eta^*_{d-i-1}(\theta^*_0)  
\tau_{h-i-1}(\theta_d) 
}
{
 \eta^*_{d-h}(\theta^*_0)  
}
p_{h+1}
\cdots
p_{d}.
\end{eqnarray*}
In the sum on the right in the above equation, we collect
terms to get a linear combination of
$\lbrace p_{h+1} \cdots p_d \rbrace_{h=i}^d$.
In this linear combination,
for $i \leq h \leq d$ let $\gamma_h$ denote the coefficient
of 
$p_{h+1} \cdots p_d$. We show
\begin{eqnarray}
\label{eq:check}
\gamma_h = 
\lbrack i,h-i,d-h\rbrack
 \frac{\eta^*_{d-i}(\theta^*_0)  
\tau_{h-i}(\theta_d) 
}
{
 \eta^*_{d-h}(\theta^*_0)  
}.
\end{eqnarray}
 First assume $h=i$. Then line
(\ref{eq:check}) holds since both sides equal $1$.
Next assume $i+1 \leq h \leq d-1$. By construction
\begin{eqnarray*}
&&\gamma_h \; =\;
\lbrack i,h-i,d-h\rbrack
a_{i+1,h+1}
 \frac{\eta^*_{d-i-1}(\theta^*_0)  
\tau_{h-i}(\theta_d) 
}
{
 \eta^*_{d-h-1}(\theta^*_0)  
} \\
&& \qquad  \quad + \quad
c_{i+1,h}
 \frac{\eta^*_{d-i-1}(\theta^*_0)  
\tau_{h-i-1}(\theta_d) 
}
{
 \eta^*_{d-h}(\theta^*_0)  
}.
\end{eqnarray*}
\noindent 
Line
(\ref{eq:check})
will follow from this
provided that
\begin{eqnarray*}
&&c_{i+1,h}
=
\lbrack i,h-i,d-h\rbrack
(\theta_d-\theta_{h-i-1})
\bigl(
\theta^*_0-\theta^*_{i+1}
-
a_{i+1,h+1}
(\theta^*_0-\theta^*_{h+1})
\bigr).
\end{eqnarray*}
The above equation is routinely checked
using
Definition \ref{lem:brackform} and
the definition of $c_{i+1,h}$ in
Lemma
\ref{lem:cij}.
We have now verified  
(\ref{eq:check}) for $i+1 \leq h \leq d-1$.
Next assume $h=d$. By construction
\begin{eqnarray*}
\gamma_d = c_{i+1,d} \eta^*_{d-i-1}(\theta^*_0)
\tau_{d-i-1}(\theta_d).
\end{eqnarray*}
Using 
the definition of
$c_{i+1,d}$ in
Lemma
\ref{lem:cij}
along with
Lemma
\ref{lem:sym}(ii) and
$p_d(\theta_0 \theta^*_0+
\theta_d \theta^*_d)=0$,
we find
$c_{i+1,d}= 
(\theta_d-\theta_{d-i-1})
(\theta^*_0-\theta^*_{i+1})$.
By these comments
\begin{eqnarray*}
\gamma_d = 
\eta^*_{d-i}(\theta^*_0)
\tau_{d-i}(\theta_d).
\end{eqnarray*}
 By the above line and Lemma
\ref{lem:sym}(ii),
line (\ref{eq:check}) holds at $h=d$.
We have now verified 
 (\ref{eq:check}) for $i \leq h \leq d$ and the result follows. 
\hfill $\Box$ \\

\begin{lemma}
\label{thm:main1}
With reference to Definitions \ref{def:aap}
and 
\ref{def:pq}, we have
 $P^\Downarrow=P$.
\end{lemma}
\noindent {\it Proof:}
By Definition
\ref{def:pq},
\begin{eqnarray}
P^\Downarrow
&=& \sum _{i=0}^d \zeta^\Downarrow_i p^\Downarrow_{i+1}
 p^\Downarrow_{i+2}\cdots p^\Downarrow_d.
\label{eq:pd}
\end{eqnarray}
In (\ref{eq:pd}), for $0 \leq i \leq d$
we expand $\zeta^\Downarrow_i$ and
$p^\Downarrow_{i+1}
 p^\Downarrow_{i+2} \cdots p^\Downarrow_d$
using Lemma
\ref{lem:zetadown}
 and
Lemma
\ref{prop:aapdvsp},
 respectively.
Now $P^\Downarrow$ becomes a linear combination of terms,
each of the form
$\zeta_r p_{s+1}p_{s+2}\cdots p_{d}$
 with $0 \leq r\leq s\leq d$.
To show $P^\Downarrow=P$, we show that for
 $0 \leq r\leq s\leq d$
the coefficient of
$\zeta_r p_{s+1}p_{s+2}\cdots p_{d}$ 
in this linear combination 
is $\delta_{r,s}$.
Using in order
(\ref{eq:pd}),
  Lemma
\ref{lem:zetadown},
Lemma 
\ref{prop:aapdvsp}, 
a change of variables $h=s-i$, 
Lemma \ref{lem:double},
Lemma \ref{lem:etatau},
Lemma \ref{lem:tauzero}(ii),
Lemma \ref{lem:sym}(ii),
the coefficient is
\begin{eqnarray*}
&&\sum_{i=r}^s  
\lbrack r,i-r,d-i\rbrack
 \frac{
\eta^*_{d-r}(\theta^*_0)  
\eta_{i-r}(\theta_0) 
}
{
\eta^*_{d-i}(\theta^*_0)
}
\lbrack i,s-i,d-s\rbrack
 \frac{
\eta^*_{d-i}(\theta^*_0)  
\tau_{s-i}(\theta_d) 
}
{
\eta^*_{d-s}(\theta^*_0)
}
\\
&& \qquad \quad = 
\quad \frac{
\eta^*_{d-r}(\theta^*_0)  
}
{
\eta^*_{d-s}(\theta^*_0)
}
\sum_{i=r}^s
\lbrack r,i-r,d-i\rbrack
\lbrack i,s-i,d-s\rbrack
\eta_{i-r}(\theta_0) 
\tau_{s-i}(\theta_d) 
\\
&& \qquad \quad = 
\quad \frac{
\eta^*_{d-r}(\theta^*_0)  
}
{
\eta^*_{d-s}(\theta^*_0)
}
\sum_{h=0}^{s-r}
\lbrack r,s-r-h,d-s+h\rbrack
\lbrack s-h,h,d-s\rbrack
\eta_{s-r-h}(\theta_0) 
\tau_{h}(\theta_d) 
\\
&& \qquad \quad = 
\quad \frac{
\eta^*_{d-r}(\theta^*_0)  
}
{
\eta^*_{d-s}(\theta^*_0)
}
\lbrack r,s-r,d-s\rbrack
\sum_{h=0}^{s-r}
\lbrack h,s-r-h,d-s+r\rbrack
\eta_{s-r-h}(\theta_0) 
\tau_{h}(\theta_d) 
\\
&& \qquad \quad = 
\quad \frac{
\eta^*_{d-r}(\theta^*_0)  
}
{
\eta^*_{d-s}(\theta^*_0)
}
\lbrack r,s-r,d-s\rbrack
\eta_{s-r}(\theta_d) 
\\
&& \qquad \quad = 
\quad 
\delta_{r,s}.
\end{eqnarray*}
The result follows.
\hfill $\Box$ \\

\begin{proposition}
\label{thm:main2}
With reference to Definitions \ref{def:aap}
and 
\ref{def:pq},
the polynomial $P$ is $D_4$-invariant.
\end{proposition}
\noindent {\it Proof:} 
This follows from 
Lemma \ref{lem:dualsame} and
Lemma \ref{thm:main1}, since
 $D_4$ is generated by
$*, \Downarrow$.
\hfill $\Box$ \\

\begin{proposition}
\label{prop:lastterm}
With reference to Definitions \ref{def:aap}
and
\ref{def:pq}, the following
 {\rm (i)}, {\rm (ii)} hold.
\begin{itemize}
\item[\rm (i)]
 $P(\theta_0 \theta^*_0+\theta_d \theta^*_d)= \zeta_d$;
\item[\rm (ii)]
 $P(\theta_0 \theta^*_d+\theta_d \theta^*_0)= \zeta^\Downarrow_d$.
\end{itemize}
\end{proposition}
\noindent {\it Proof:}
(i)
Assume $d \geq 1$; otherwise
both sides equal $1$.
Now evaluate
(\ref{eq:drinp})
 at $\lambda=\theta_0 \theta^*_0 + \theta_d \theta^*_d$
 and observe 
$p_d(
 \theta_0 \theta^*_0 + \theta_d \theta^*_d
)=0$ by Lemma
\ref{lem:dnorm}.
\\
\noindent (ii) Apply (i) to $\Phi^\Downarrow$.
\hfill $\Box$ \\

\section{The Drinfel'd polynomial for type ${\rm III}^+$}

\noindent 
In this section we define the Drinfel'd polynomial
for a sharp tridiagonal system of type 
${\rm III}^+$. In order to motivate our 
definition, consider the type 
 I expression for $\alpha_i$ that appears in the
table of Definition 
\ref{def:aap}. That expression is a fraction
with $(q^i-q^{-i})^2$ in the numerator
and 
$(q^d-q^{-d})^2$ in the denominator.
Recall that $q^2+q^{-2}=\beta$ where $\beta$ is the type of
$\Phi$.
Under the  assumption of type 
 ${\rm III}^+$ we have $\beta = -2$
so $q^2$ becomes $-1$.
In this case the 
expression $(q^i-q^{-i})^2$
 becomes $0$ if $i$ is even
and $-4$ if $i$ is odd.
Since $d$ is even
the expression $(q^d-q^{-d})^2$ becomes $0$.
For $i$ even we can ``take limits'' 
and evaluate $\alpha_i$ in a reasonable way, but 
if we try this for $i$ odd then $\alpha_i$ becomes $\infty$.
To resolve this problem
we multiply the polynomial
$P$ in Definition \ref{def:pq} by $(q^d-q^{-d})^d$
before taking limits.
This procedure gives a reasonable limit
but it does have the effect of sending the $i$-summand in
(\ref{eq:drinp}) to zero for $1 \leq i \leq d$.
Consequently for type 
${\rm III}^+$ the Drinfel'd polynomial is no longer defined as a sum.
Instead it looks as follows.

\begin{definition}
\label{def:3p}
\rm
Let $\Phi$ denote a sharp tridiagonal system over 
$\K$ with eigenvalue sequence
 $\lbrace \theta_i\rbrace_{i=0}^d$
and dual eigenvalue sequence  
 $\lbrace \theta^*_i\rbrace_{i=0}^d$.
Assume $\Phi$ 
 has type ${\rm III}^{+}$.
 We define
the Drinfel'd polynomial of $\Phi$ to be
\begin{eqnarray}
P = p_1 p_2 \cdots p_d
\label{eq:p3}
\end{eqnarray}
where
for $1 \leq i \leq d$,
\begin{eqnarray}
\label{eq:case3p}
p_i &=&     
\begin{cases}  
(\theta_0 \theta^*_d + \theta_d \theta^*_0 - \lambda)i^2/d^2+
(\theta_0-\theta_i)(\theta^*_0-\theta^*_i) 
& \text{\rm if $i$ is even},
\\
\theta_0 \theta^*_d + \theta_d \theta^*_0 - \lambda 
& \text{\rm if $i$ is odd}.
    \end{cases} 
\end{eqnarray}
Note that $P=1$ if $d=0$.
\end{definition}

\noindent Referring to Definition
\ref{def:3p}, we now put $P$ in a more attractive form.

\begin{lemma}
\label{lem:pialt}
With reference to Definition
\ref{def:3p},
 for even $i$ $(1 \leq i \leq d)$ we have
\begin{eqnarray*} 
p_i= 
(\theta_0 \theta^*_0 + \theta_d \theta^*_d -\lambda)i^2/d^2.
\end{eqnarray*}
\end{lemma} 
\noindent {\it Proof:}
By Lemma \ref{lem:oioiIIIp},
\begin{eqnarray*}
(\theta_0-\theta_i)
(\theta^*_0-\theta^*_i) = 
(\theta_0-\theta_d)
(\theta^*_0-\theta^*_d)i^2/d^2.
\end{eqnarray*}
Evaluating the first line of 
(\ref{eq:case3p}) using this we get the result.
\hfill $\Box$ \\

\begin{lemma}
\label{lem:ptype3p}
With reference to Definition
\ref{def:3p}
and in the notation of
 Lemma \ref{lem:eig3}, for $1 \leq i \leq d$ we have
\begin{eqnarray*}
p_i &=&     
\begin{cases}  
(2(a+b)(a^*+b^*)+cc^*d^2/2 - \lambda)i^2/d^2
& \text{\rm if $i$ is even},
\\
2(a+b)(a^*+b^*)-cc^*d^2/2 - \lambda
& \text{\rm if $i$ is odd}.
    \end{cases} 
\end{eqnarray*}
\end{lemma}
\noindent {\it Proof:}
Evaluate 
(\ref{eq:case3p})
using Lemma
 \ref{lem:eig3}
and 
Lemma \ref{lem:pialt}.
\hfill $\Box$ \\

\begin{lemma}
\label{prop:type3p}
With reference to Definition \ref{def:3p},
assume $d\geq 2$ and abbreviate $D=d/2$.
\begin{itemize}
\item[{\rm (i)}]
 $P$ looks
as follows
in terms of the eigenvalues and dual eigenvalues:
\begin{eqnarray*}
P \; =\;
D^{-d}\,(D!)^2\,
 (\theta_0 \theta^*_0 + \theta_d \theta^*_d- \lambda)^{d/2} 
\; (\theta_0 \theta^*_d + \theta_d \theta^*_0- \lambda)^{d/2}.
\end{eqnarray*}
\item[{\rm (ii)}]
$P$ looks as follows in the 
notation of
 Lemma \ref{lem:eig3}:
\begin{eqnarray*}
&& P \; =\;
D^{-d}\,(D!)^2\,
 \bigl(
2(a+b)(a^*+b^*)+cc^*d^2/2 - \lambda \bigr)^{d/2}
\\
&& \qquad \qquad \quad \qquad \qquad \qquad 
\times \quad 
\bigl(2(a+b)(a^*+b^*)-cc^*d^2/2 - \lambda \bigr)^{d/2}.
\end{eqnarray*}
\end{itemize}
\end{lemma}
\noindent {\it Proof:}
(i) 
 Evaluate
(\ref{eq:p3}) using
(\ref{eq:case3p}) 
and
Lemma
\ref{lem:pialt}.
\\
(ii) Evaluate 
(\ref{eq:p3}) using
Lemma \ref{lem:ptype3p}.
\hfill $\Box$ \\

\begin{proposition}
\label{prop:3d4}
With reference to Definition \ref{def:3p},
$P$ is $D_4$-invariant.
\end{proposition}
\noindent {\it Proof:}
Use Lemma \ref{prop:type3p}(i).
\hfill $\Box$ \\

\begin{proposition}
With reference to Definition \ref{def:3p}, 
the following
 {\rm (i)}, {\rm (ii)} hold for $d\geq 2$.
\begin{itemize}
\item[\rm (i)]
 $P(\theta_0 \theta^*_0+\theta_d \theta^*_d)= 0$;
\item[\rm (ii)]
 $P(\theta_0 \theta^*_d+\theta_d \theta^*_0)= 0$.
\end{itemize}
\end{proposition}
\noindent {\it Proof:}
Immediate from Lemma 
\ref{prop:type3p}(i).
\hfill $\Box$ \\

\noindent Combining Proposition
\ref{thm:main2}
and Proposition 
\ref{prop:3d4}
we obtain the following theorem.

\begin{theorem}
\label{thm:maind4}
Let $\Phi$ denote a sharp tridiagonal system
and let $P$ denote the corresponding Drinfel'd polynomial
from Definition 
\ref{def:pq}
and Definition \ref{def:3p}.
Then $P$ is $D_4$-invariant.
\end{theorem}

\begin{definition}
\label{def:drintdpair}
\rm
Let $A,A^*$ denote a sharp tridiagonal pair.
By the {\it Drinfel'd polynomial} of $A,A^*$
we mean the Drinfel'd polynomial of an
associated tridiagonal system.
By construction $A,A^*$ and $A^*,A$ have
the same Drinfel'd polynomial.
\end{definition}

\section{The normalized Drinfel'd polynomial for type I}

\noindent 
Let $\Phi$ denote a sharp tridiagonal system over $\K$.
In Definition 
\ref{def:pq} and
Definition
 \ref{def:3p}
 we defined the Drinfel'd polynomial
$P$ of $\Phi$. 
One advantage of our definition is that it depends in
only a minor way on the type of $\Phi$.
 One disadvantage is that the roots
of $P$ are not as nice as they could be.
We make the roots nicer by 
 introducing the 
{\it normalized Drinfel'd polynomial} $\hat P$ on a type-by-type basis.
For each type,  $\hat P$ is related to $P$
 by an equation
${\hat P}(\lambda)= P(u\lambda+v)$,
where $u,v$ are scalars in $\overline \K$ 
that depend on the type. For the rest of
 this section we focus on type I.

\begin{assumption}
\label{def:a1}
\rm
 Let $\Phi$ denote a sharp tridiagonal system
over $\K$, with eigenvalue sequence
$\lbrace \theta_i\rbrace_{i=0}^d$ and
dual eigenvalue sequence
$\lbrace \theta^*_i\rbrace_{i=0}^d$. 
Assume $\Phi$ has type I and let the scalars
$q,a,b,c,a^*,b^*,c^*$ be as in
Lemma
 \ref{lem:eig1}.
\end{assumption}

\medskip
\noindent The following definition is
motivated by 
Lemma
\ref{lem:oioiI}.

\begin{definition}
\label{def:pialt}
\rm
With reference to Assumption \ref{def:a1},
for $1 \leq i \leq d$ 
define 
${\hat p}_i \in \K \lbrack \lambda \rbrack$ by
\begin{eqnarray}
{\hat p}_i = (q^i-q^{-i})^2(  
bb^*q^{2i-2d}+cc^*q^{2d-2i}-\lambda).
\label{eq:pinormcase1}
\end{eqnarray}
By Lemma \ref{lem:oioiI},
\begin{eqnarray}
{\hat p}_i(bc^*+cb^*)= (\theta_0-\theta_i)(\theta^*_0-\theta^*_i).
\label{eq:hatpnorm}
\end{eqnarray} 
\end{definition}

\begin{definition}
\label{def:palt}
\rm
With reference to Assumption \ref{def:a1},
 define 
${\hat P} \in \K \lbrack \lambda \rbrack$ by
\begin{eqnarray}
\label{eq:drinnorm}
{\hat P} = \sum _{i=0}^d \zeta_i {\hat p}_{i+1}
 {\hat p}_{i+2} \cdots {\hat p}_d,
\end{eqnarray}
where $\lbrace \zeta_i\rbrace_{i=0}^d$ is the split
sequence of $\Phi$ and the ${\hat p}_i$ are from
Definition
\ref{def:pialt}. We call
$\hat P$ the {\it normalized Drinfel'd polynomial} for $\Phi$.
\end{definition}

\begin{lemma}
\label{lem:rs}
With reference to Assumption \ref{def:a1}, both
\begin{eqnarray*}
u (b b^* + c c^*)+ v &=& \theta_0  \theta^*_0 + 
\theta_d  \theta^*_d, \\
u (b c^* + c b^*)+ v &=& \theta_0  \theta^*_d + 
\theta_d  \theta^*_0,
\end{eqnarray*}
where
\begin{eqnarray*}
u &=& (q^d-q^{-d})^2,
\label{eq:rdef}
\\
v &=& 2a a^*+ 2(b+c)(b^*+c^*)  + (q^d+q^{-d})a(b^*+c^*)
 + (q^d+q^{-d})a^*(b+c).
\label{eq:sdef}
\end{eqnarray*} 
\end{lemma}
\noindent {\it Proof:} Use
Lemma
 \ref{lem:eig1}.
\hfill $\Box$ \\

\begin{theorem}
\label{thm:hatp}
With reference to Assumption \ref{def:a1},
\begin{eqnarray}
{\hat p}_i(\lambda) &=& p_i(u \lambda + v) \qquad \qquad (1 \leq i \leq d),
\label{eq:pihat}
\\
{\hat P}(\lambda) &=& P(u \lambda + v),
\label{eq:phat}
\end{eqnarray}
where $u,v$ are from Lemma
\ref{lem:rs}.
\end{theorem}
\noindent {\it Proof:}
Each side of
(\ref{eq:pihat})
is a polynomial in 
$\K\lbrack \lambda \rbrack$ with degree 1.
On each side of
(\ref{eq:pihat})
the polynomial has $\lambda $ coefficient  $-(q^i-q^{-i})^2$.
Moreover at
 $\lambda=bc^*+cb^*$ this polynomial takes the
value
$(\theta_0-\theta_i)(\theta^*_0-\theta^*_i)$.
We have verified (\ref{eq:pihat}) and
 (\ref{eq:phat}) follows.
\hfill $\Box$ \\

\begin{proposition}
\label{prop:d4I}
With reference to Assumption \ref{def:a1} and
Definition \ref{def:palt}, the polynomial 
 $\hat P$ is $D_4$-invariant.
\end{proposition}
\noindent {\it Proof:}
Follows from
Proposition
\ref{thm:main2}
and (\ref{eq:phat}).
\hfill $\Box$ \\

\begin{proposition}
\label{prop:lasttermI}
With reference to Assumption \ref{def:a1}
and Definition
\ref{def:palt}, the following
{\rm (i)}, {\rm (ii)} hold.
\begin{itemize}
\item[\rm (i)]
 ${\hat P}(b b^* + c c^*)= \zeta_d$;
\item[\rm (ii)]
 ${\hat P}(b c^* + c b^*)= \zeta^\Downarrow_d$.
\end{itemize}
\end{proposition}
\noindent {\it Proof:}
Combine
Proposition
\ref{prop:lastterm},
Lemma \ref{lem:rs},
and 
(\ref{eq:phat}).
\hfill $\Box$ \\

\noindent We now consider the normalized Drinfel'd polynomial
for a Leonard system 
of type I.

\begin{lemma}
\label{ex:lpI}
{\rm \cite[Theorem~6.1]{NT:balanced}}
With reference to Assumption \ref{def:a1},
suppose $\Phi$ is a Leonard system. Let
$\lbrace \varphi_i\rbrace_{i=1}^d$
(resp. 
$\lbrace \phi_i\rbrace_{i=1}^d$) denote the corresponding
first (resp. second) split sequence, from
 \cite[Section~3]{LS99}. 
Then there exists $t \in {\overline \K}$ such that
for $1 \leq i \leq d$,
\begin{eqnarray*}
\varphi_i &=& (q^i-q^{-i}) 
(q^{d-i+1}-q^{i-d-1})(t-b b^* q^{2i-d-1}- 
c c^* q^{d+1-2i}),
\\
\phi_i &=& (q^i-q^{-i}) 
(q^{d-i+1}-q^{i-d-1})(t-c b^* q^{2i-d-1}- 
b c^* q^{d+1-2i}).
\end{eqnarray*}
The above expressions can be factored more completely if
$b b^* c c^*\not=0$. In this case
\begin{eqnarray*}
\varphi_i &=& (q^i-q^{-i}) 
(q^{d-i+1}-q^{i-d-1})(q^{-i}-b b^* \psi^{-1} q^{i-d-1}) 
(q^i \psi-c c^*  q^{d-i+1}),
\\
\phi_i &=& (q^i-q^{-i}) 
(q^{d-i+1}-q^{i-d-1})(q^{-i}-c b^* \psi^{-1} q^{i-d-1}) 
(q^{i} \psi-b c^* q^{d-i+1}),
\end{eqnarray*}
where $\psi \in {\overline \K}$ is a solution to
\begin{eqnarray*}
\psi + b b^* c c^* \psi^{-1} = t.
\end{eqnarray*}
\end{lemma}

\begin{proposition}
\label{prop:rootsI}
With reference to Assumption \ref{def:a1},
suppose $\Phi$ is a Leonard system and
let $\hat P$ denote the corresponding normalized Drinfel'd
polynomial.
\begin{itemize}
\item[{\rm (i)}] Assume $b b^* c c^* \not=0$. Then the roots of
$\hat P$ are 
\begin{eqnarray*}
\psi q^{d+1-2i} + b b^* c c^* \psi^{-1} q^{2i-d-1}    
\qquad \qquad (1 \leq i \leq d),
\end{eqnarray*}
where $\psi $ is from Lemma \ref{ex:lpI}.
\item[{\rm (ii)}] Assume $b b^* c c^* = 0$. Then the roots of $\hat P$
are 
\begin{eqnarray*}
t q^{d+1-2i} \qquad \qquad (1 \leq i \leq d),
\end{eqnarray*}
where $t$ is from Lemma
\ref{ex:lpI}.
\end{itemize}
\end{proposition}
\noindent {\it Proof:}
Assume $d\geq 1$ to avoid trivialities.
Let 
$\lbrace \zeta_i\rbrace_{i=0}^d$ denote the split
sequence of $\Phi$ from Section 6,
and let $\lbrace \varphi_i\rbrace_{i=1}^d$
denote the first split sequence of $\Phi$.
By Note
\ref{nt:spltdiff},
$\zeta_i = \varphi_1\varphi_2 \cdots \varphi_i$ 
for $0 \leq i \leq d$.
By this and
 (\ref{eq:drinnorm}),
$\hat P$ is equal to 
$\varphi_1 \varphi_2 \cdots \varphi_d $ times
\begin{eqnarray}
\label{eq:adj}
\sum_{n=0}^d \frac{
{\hat p}_d {\hat p}_{d-1} \cdots {\hat p}_{d-n+1}
}
{\varphi_d \varphi_{d-1} \cdots \varphi_{d-n+1}}.
\end{eqnarray}
The 
denominators in (\ref{eq:adj}) are nonzero, since
each of $\varphi_1, \varphi_2, \ldots, \varphi_d$ is nonzero
by 
\cite[Theorem~1.9]{LS99}.
\\
\noindent (i)
We evaluate 
(\ref{eq:adj}) using
(\ref{eq:pinormcase1})
and 
the data in Lemma 
\ref{ex:lpI}. 
This calculation shows that (\ref{eq:adj}) is equal to
\begin{eqnarray*}
\sum_{n=0}^d \frac{(q^{-2d};q^2)_n (x/bb^*;q^2)_n 
(cc^*/x;q^2)_n q^{2n}} 
{(q^{1-d}\psi/bb^*;q^2)_n (q^{1-d}cc^*/ \psi;q^2)_n (q^2;q^2)_n},
\label{eq:phyper}
\end{eqnarray*}
where $\lambda = x + bb^*cc^*x^{-1}$.
Basic hypergeometric series are defined in \cite[p.~4]{GR}.
By that definition
the above sum
is the basic hypergeometric series
\begin{eqnarray}
 {}_3\phi_2 \biggl\lbrack{{q^{-2d}, \;x/bb^*,\;cc^*/x}\atop
{q^{1-d}\psi/bb^*,\;\;q^{1-d}cc^*/\psi}}\;  ; \; q^2,\;q^2 \biggr\rbrack.
\label{eq:3p2}
\end{eqnarray}
By the $q$-Saalsch\"{u}tz formula
\cite[p.~355]{GR} the sum
(\ref{eq:3p2}) is equal to
\begin{eqnarray}
\frac{
(q^{1-d}x/\psi;q^2)_d 
(q^{1-d}bb^*cc^*/\psi x;q^2)_d 
}
{
(q^{1-d}bb^*/\psi;q^2)_d (q^{1-d}cc/\psi;q^2)_d
}.
\label{eq:ssform}
\end{eqnarray}
By (\ref{eq:aqn}) the numerator in
(\ref{eq:ssform}) 
is equal to
\begin{eqnarray}
\label{eq:numpre}
\prod_{i=1}^d (1-q^{2i-d-1}x/ \psi)(1-q^{2i-d-1}bb^*cc^*/\psi x).
\end{eqnarray}
For $1 \leq i \leq d$ the $i$-factor in
(\ref{eq:numpre}) is equal to $\psi^{-1}q^{2i-d-1}$ times
\begin{eqnarray*}
\psi q^{d+1-2i}+bb^*cc^* \psi^{-1}q^{2i-d-1}-\lambda.
\end{eqnarray*}
Therefore the numerator in (\ref{eq:ssform}) 
is a nonzero scalar multiple
of 
\begin{eqnarray*}
\prod_{i=1}^d 
(\psi q^{d+1-2i}+bb^*cc^* \psi^{-1}q^{2i-d-1}-\lambda).
\end{eqnarray*}
The result follows.
\\
\noindent (ii) Replacing $q$ by $q^{-1}$ if necessary,
we may assume without loss that $bb^*=0$. For the moment further 
assume that the scalar $t$ from Lemma
\ref{ex:lpI} is nonzero.
Proceeding as in (i) above,
we find that (\ref{eq:adj}) is equal to
\begin{eqnarray*}
 {}_2\phi_1 \biggl\lbrack{{q^{-2d}, \;cc^*/\lambda }\atop
{q^{1-d}cc^*/t}}\;  ; \; q^2,\;\frac{q^{d+1}\lambda}{t} \biggr\rbrack,
\end{eqnarray*}
which is equal to 
\begin{eqnarray}
\frac{
(q^{1-d}\lambda/t;q^2)_d 
}
{
(q^{1-d}cc^*/t;q^2)_d
}
\label{eq:ssform2}
\end{eqnarray}
by the $q$-Chu-Vandermonde formula
\cite[p.~354]{GR}.
By (\ref{eq:aqn}) the numerator in
(\ref{eq:ssform2}) 
is a nonzero scalar multiple
of 
\begin{eqnarray*}
\prod_{i=1}^d 
(t q^{d+1-2i}-\lambda),
\end{eqnarray*}
giving the result for the case $t\not=0$.
Next assume $t=0$. Then $cc^*\not=0$; otherwise
$\varphi_1=0$.
In this case  (\ref{eq:adj}) is equal to
\begin{eqnarray*}
 {}_2\phi_1 \biggl\lbrack{{q^{2d}, \;\lambda/cc^*}\atop
{0\;}} \; ; \;q^{-2},\;q^{-2} \biggr\rbrack,
\end{eqnarray*}
which is equal to 
$(\lambda/cc^*)^d$
by another version of the  $q$-Chu-Vandermonde formula
\cite[p.~354]{GR}.
The result follows for the case $t=0$, 
 and the
proof is complete.
\hfill $\Box$ \\

\section{The normalized Drinfel'd polynomial for type II}

In this section we introduce
the normalized Drinfel'd polynomial for a sharp
tridiagonal system of type II. 

\begin{assumption}
\label{def:a2}
\rm
 Let $\Phi$ denote a sharp tridiagonal system
over $\K$, with eigenvalue sequence
$\lbrace \theta_i\rbrace_{i=0}^d$ and
dual eigenvalue sequence
$\lbrace \theta^*_i\rbrace_{i=0}^d$. 
Assume $\Phi$ has type II and let the scalars
$a,b,c,a^*,b^*,c^*$ be as in
Lemma
\ref{lem:eig2}.
\end{assumption}

\medskip
\noindent The following definition is motivated by 
Lemma
\ref{lem:oioiII}.

\begin{definition}
\label{def:pialtII}
\rm
With reference to Assumption \ref{def:a2},
for $1 \leq i \leq d$ 
 define 
${\hat p}_i \in \K \lbrack \lambda \rbrack$ by
\begin{eqnarray}
\label{eq:normdp2}
{\hat p}_i = i^2(  
 b b^*/2 + (b c^* + c b^*)(d-i) + c c^* (d-i)^2 -\lambda ).
\end{eqnarray}
By Lemma
\ref{lem:oioiII},
\begin{eqnarray}
{\hat p}_i(-bb^*/2)= (\theta_0-\theta_i)(\theta^*_0-\theta^*_i).
\label{eq:hatpnormII}
\end{eqnarray} 
\end{definition}

\begin{definition}
\label{def:paltII}
\rm
With reference to Assumption \ref{def:a2},
 define
${\hat P} \in \K \lbrack \lambda \rbrack$ by
\begin{eqnarray}
\label{eq:drinnormII}
{\hat P} = \sum _{i=0}^d \zeta_i {\hat p}_{i+1}
 {\hat p}_{i+2} \cdots {\hat p}_d,
\end{eqnarray}
where $\lbrace \zeta_i\rbrace_{i=0}^d$ is the split
sequence of $\Phi$ and the ${\hat p}_i$ are from
Definition
\ref{def:pialtII}. We call
$\hat P$ the {\it normalized Drinfel'd polynomial} for $\Phi$.
\end{definition}

\begin{lemma}
\label{lem:rsII}
With reference to Assumption \ref{def:a2}, both
\begin{eqnarray*}
u b b^*/2+ v &=& \theta_0  \theta^*_0 + 
\theta_d  \theta^*_d, \\
-u b b^*/2+ v &=& \theta_0  \theta^*_d + 
\theta_d  \theta^*_0,
\end{eqnarray*}
where 
\begin{eqnarray*}
u &=& d^2,
\label{eq:rdefII}
\\
v &=& 2a a^*.
\label{eq:sdefII}
\end{eqnarray*} 
\end{lemma}
\noindent {\it Proof:}
Use Lemma
\ref{lem:eig2}.
\hfill $\Box$ \\

\begin{theorem}
\label{thm:hatpII}
With reference to Assumption \ref{def:a2},
\begin{eqnarray}
{\hat p}_i(\lambda) &=& p_i(u \lambda + v) \qquad \qquad (1 \leq i \leq d),
\label{eq:pihatII}
\\
{\hat P}(\lambda) &=& P(u \lambda + v),
\label{eq:phatII}
\end{eqnarray}
where $u,v$ are from Lemma
\ref{lem:rsII}.
\end{theorem}
\noindent {\it Proof:}
Each side of 
(\ref{eq:pihatII}) is a polynomial in 
$\K\lbrack \lambda \rbrack$ with degree 1.
On each side of
(\ref{eq:pihatII})
the polynomial has $\lambda $ coefficient  $-i^2$.
Moreover at
 $\lambda=-bb^*/2$ this polynomial takes the
value
$(\theta_0-\theta_i)(\theta^*_0-\theta^*_i)$.
We have verified (\ref{eq:pihatII}) and
(\ref{eq:phatII}) follows.
\hfill $\Box$ \\

\begin{proposition}
\label{prop:d4II}
With reference to Assumption \ref{def:a2}
and Definition
\ref{def:paltII}, the polynomial $\hat P$
is $D_4$-invariant.
\end{proposition}
\noindent {\it Proof:}
Follows from
Proposition
\ref{thm:main2}
and (\ref{eq:phatII}).
\hfill $\Box$ \\

\begin{proposition}
\label{prop:lasttermII}
With reference to Assumption \ref{def:a2}
and Definition
\ref{def:paltII},
the following
 {\rm (i)}, {\rm (ii)} hold.
\begin{itemize}
\item[\rm (i)]
 ${\hat P}(b b^*/2)= \zeta_d$;
\item[\rm (ii)]
 ${\hat P}(-b b^*/2)= \zeta^\Downarrow_d$.
\end{itemize}
\end{proposition}
\noindent {\it Proof:}
Combine
Proposition
\ref{prop:lastterm},
Lemma \ref{lem:rsII},
and 
(\ref{eq:phatII}).
\hfill $\Box$ \\

\noindent We now consider the normalized Drinfel'd polynomial
for a Leonard system of type II.

\begin{lemma}
\label{ex:lpII}
{\rm \cite[Theorem~7.1]{NT:balanced}}
With reference to Assumption \ref{def:a2},
suppose $\Phi$ is a Leonard system. Let
$\lbrace \varphi_i\rbrace_{i=1}^d$
(resp. 
$\lbrace \phi_i\rbrace_{i=1}^d$) denote the corresponding
first (resp. second) split sequence, from
 \cite[Section~3]{LS99}. 
Then there exists $t \in {\overline \K}$ such that
for $1 \leq i \leq d$,
\begin{eqnarray*}
\varphi_i &=& i 
(d-i+1)\bigl(t-b b^*/2 
+(b c^* + c b^*)(i-\tfrac{d+1}{2})- c c^* (i-\tfrac{d+1}{2})^2\bigr),
\\
\phi_i &=&
i (d-i+1)\bigl(t+b b^*/2 
+(c b^*-b c^*)(i-\tfrac{d+1}{2})- c c^* (i-\tfrac{d+1}{2})^2\bigr).
\end{eqnarray*}
The above expressions can be factored more completely if
$c c^*\not=0$. In this case
\begin{eqnarray*}
\varphi_i &=&  
c^{-1} c^{*-1} i (d-i+1)
\bigl(\tfrac{\psi + b c^* + c b^*}{2} - c c^* (i-\tfrac{d+1}{2})\bigr)
\bigl(\tfrac{\psi -b c^* - c b^*}{2} + c c^* (i-\tfrac{d+1}{2})\bigr),
\\
\phi_i &=&  
c^{-1} c^{*-1} i (d-i+1)
\bigl(\tfrac{\psi - b c^*+ c b^*}{2} - c c^* (i-\tfrac{d+1}{2})\bigr)
\bigl(\tfrac{\psi + b c^* -c b^*}{2} + c c^* (i-\tfrac{d+1}{2})\bigr),
\end{eqnarray*}
where $\psi \in {\overline \K}$ is a solution to
\begin{eqnarray*}
\psi^2 = 4 t c c^* + b^2 c^{*2} + b^{* 2} c^2.
\end{eqnarray*}
\end{lemma} 

\begin{proposition}
\label{prop:rootsII}
With reference to Assumption \ref{def:a2},
suppose $\Phi$ is a Leonard system  and
let $\hat P$ denote the corresponding normalized Drinfel'd
polynomial.
\begin{itemize}
\item[{\rm (i)}] Assume $ c c^* \not=0$. Then the roots of
$\hat P$ are 
\begin{eqnarray*}
t + \psi (i- \tfrac{d+1}{2})+ c c^* (i- \tfrac{d+1}{2})^2
\qquad \qquad (1 \leq i \leq d),
\end{eqnarray*}
where $\psi, t $ are from Lemma \ref{ex:lpII}.
\item[{\rm (ii)}] Assume $ c c^* = 0$. Then the roots of $\hat P$
are 
\begin{eqnarray*}
t + (b c^*+ c b^*)(i-\tfrac{d+1}{2}) 
\qquad \qquad (1 \leq i \leq d),
\end{eqnarray*}
where $t$ is from Lemma
\ref{ex:lpII}.
\end{itemize}
\end{proposition}
\noindent {\it Proof:}
We begin as in the proof of
Proposition
\ref{prop:rootsI}.
\\
\noindent (i)
We evaluate 
(\ref{eq:adj}) using
(\ref{eq:normdp2}) 
and 
the data in Lemma 
\ref{ex:lpII}. 
The result is that (\ref{eq:adj}) is equal to
\begin{eqnarray*}
\sum_{n=0}^d \frac{
(-d)_n
\bigl(\tfrac{bc^*+cb^*+x}{2cc^*} \bigr)_n
\bigl(\tfrac{bc^*+cb^*-x}{2cc^*}\bigr)_n
 } 
{
\bigl(\tfrac{bc^*+cb^*+ \psi}{2cc^*}+ \tfrac{1-d}{2}\bigr)_n
\bigl(\tfrac{bc^*+cb^*-\psi}{2cc^*}+\tfrac{1-d}{2} \bigr)_n 
n!
},
\label{eq:Fhyper}
\end{eqnarray*}
where $x^2= 4\lambda cc^*+ b^2c^{*2}+ b^{*2}c^2$
and 
$a_n=a(a+1)\cdots (a+n-1)$.
Hypergeometric series are defined in \cite[p.~3]{GR}.
By that definition
the above sum is the hypergeometric series
\begin{eqnarray}
 {}_3F_2 \Biggl\lbrack{{-d,\;\; \;
\tfrac{bc^*+cb^*+x}{2cc^*}, \;\; \;
\tfrac{bc^*+cb^*-x}{2cc^*}}\atop
{\tfrac{bc^*+cb^*+\psi}{2cc^*}+\frac{1-d}{2},\;\;
\tfrac{bc^*+cb^*-\psi}{2cc^*}+\tfrac{1-d}{2}} 
}\;  ; \; 1 \Biggr\rbrack.
\label{eq:3F2}
\end{eqnarray}
By the Saalsch\"{u}tz formula
\cite[p.~17]{GR} the sum
(\ref{eq:3F2}) is equal to
\begin{eqnarray}
\frac{
\bigl(\tfrac{\psi-x}{2cc^*}+\tfrac{1-d}{2}\bigr)_d 
\bigl(\tfrac{\psi+x}{2cc^*}+\tfrac{1-d}{2}\bigr)_d 
}
{
\bigl(\tfrac{\psi+bc^*+cb^*}{2cc^*}+\tfrac{1-d}{2}\bigr)_d
\bigl(\tfrac{\psi-bc^*-cb^*}{2cc^*}+\tfrac{1-d}{2}\bigr)_d
}.
\label{eq:ssformII}
\end{eqnarray}
The numerator in
(\ref{eq:ssformII}) 
is equal to
\begin{eqnarray}
\label{eq:numpreII}
\prod_{i=1}^d 
\bigl(\tfrac{\psi-x}{2cc^*}+i-\tfrac{d+1}{2}\bigr)
\bigl(\tfrac{\psi+x}{2cc^*}+i-\tfrac{d+1}{2}\bigr). 
\end{eqnarray}
For $1 \leq i \leq d$ the $i$-factor in
(\ref{eq:numpreII}) is equal to $(cc^*)^{-1}$ times
\begin{eqnarray*}
 t +\psi (i- \tfrac{d+1}{2})+ c c^* (i- \tfrac{d+1}{2})^2-\lambda.
\end{eqnarray*}
Therefore the numerator in (\ref{eq:ssformII}) 
is a nonzero scalar multiple
of 
\begin{eqnarray*}
\prod_{i=1}^d 
\bigl(t +\psi (i- \tfrac{d+1}{2})+ c c^* (i- \tfrac{d+1}{2})^2-\lambda \bigr).
\end{eqnarray*}
The result follows.
\\
\noindent (ii) Replacing $\Phi$ by $\Phi^*$ if necessary,
 we may assume without loss that $c=0$.
 For the moment
assume further that $c^*\not=0$. Then 
 (\ref{eq:adj}) is equal to
\begin{eqnarray*}
 {}_2F_1 \Biggl\lbrack{{-d, \;\; \;
\tfrac{b^*}{2c^*}-\tfrac{\lambda}{bc^*}
}\atop
{\tfrac{b^*}{2c^*}-\tfrac{t}{bc^*}+\tfrac{1-d}{2}
}}\;  ; \; 1 \Biggr\rbrack,
\label{eq:2f1}
\end{eqnarray*}
which is equal to 
\begin{eqnarray}
\frac{
\bigl(\tfrac{\lambda-t}{bc^*}+\tfrac{1-d}{2}\bigr)_d
}
{\bigl(\tfrac{b^*}{2c^*}-\tfrac{t}{bc^*}+\tfrac{1-d}{2}\bigr)_d
}
\label{eq:ssform2II}
\end{eqnarray}
by the Chu-Vandermonde formula
\cite[p.~2]{GR}.
The numerator in
(\ref{eq:ssform2II}) 
is a nonzero scalar multiple
of 
\begin{eqnarray*}
\prod_{i=1}^d 
\bigl(t + bc^*(i-\tfrac{d+1}{2})-\lambda\bigr),
\end{eqnarray*}
giving the result for the case $c^*\not=0$.
Next assume $c^*=0$. Then $bb^*\not=2t$; otherwise
$\varphi_1=0$.
In this case  (\ref{eq:adj}) is equal to
\begin{eqnarray*}
 {}_1F_0 \biggl\lbrack{{-d}\atop
{ }} \; ; \;\frac{bb^*-2\lambda}{bb^*-2t} \biggr\rbrack,
\end{eqnarray*}
which is equal to 
$2^d(\lambda-t)^d(bb^*-2t)^{-d}$
by the binomial theorem.
The result follows for the case $c^*=0$, 
 and the
proof is complete.
\hfill $\Box$ \\

\section{The normalized Drinfel'd polynomial for type
${\rm III}$ }

\noindent In this section we introduce
the normalized Drinfel'd polynomial for a sharp
tridiagonal system of type III. 

\medskip
\noindent For a sharp tridiagonal system of type ${\rm III}^{+}$
we define the normalized Drinfel'd polynomial $\hat P$ to be the
Drinfel'd polynomial $P$ from
Definition 
\ref{def:3p}.
For the rest of this section we focus on type 
 ${\rm III}^{-}$.

\begin{assumption}
\label{def:a3}
\rm
Let $\Phi$ denote a sharp tridiagonal system
over $\K$, with eigenvalue sequence
$\lbrace \theta_i\rbrace_{i=0}^d$ and
dual eigenvalue sequence
$\lbrace \theta^*_i\rbrace_{i=0}^d$. 
Assume $\Phi$ has type ${\rm III}^{-}$
and let the scalars
$a,b,c,a^*,b^*,c^*$ be as in
Lemma
\ref{lem:eig3}.
\end{assumption}

\noindent The following definition is
motivated by 
Lemma
\ref{lem:oioiIIIp}.

\begin{definition}
\label{def:pialtIII}
\rm
With reference to Assumption \ref{def:a3},
for $1 \leq i \leq d$ 
 define 
${\hat p}_i \in \K \lbrack \lambda \rbrack$ by
\begin{eqnarray}
\label{eq:pinormd3}
{\hat p}_i = 
\begin{cases}  
cc^*i^2
 & \text{\rm if $i$ is even},  \\
2 b b^* +2(b c^*+ c b^*)(i-d)+cc^*(i-d)^2 - \lambda 
& \text{\rm if $i$ is odd},
    \end{cases} 
\end{eqnarray}
By Lemma
\ref{lem:oioiIIIp},
\begin{eqnarray}
{\hat p}_i(-2bb^*)= (\theta_0-\theta_i)(\theta^*_0-\theta^*_i).
\label{eq:hatpnormIII}
\end{eqnarray} 
\end{definition}

\begin{definition}
\label{def:paltIII}
\rm
With reference to Assumption \ref{def:a3},
 define
${\hat P} \in \K \lbrack \lambda \rbrack$ by
\begin{eqnarray}
\label{eq:drinnormIII}
{\hat P} = \sum _{i=0}^d \zeta_i {\hat p}_{i+1}
 {\hat p}_{i+2} \cdots {\hat p}_d,
\end{eqnarray}
where $\lbrace \zeta_i\rbrace_{i=0}^d$ is the split
sequence of $\Phi$ and the ${\hat p}_i$ are from
Definition
\ref{def:pialtIII}. We call
$\hat P$ the {\it normalized Drinfel'd polynomial} for $\Phi$.
We note that $\hat P$ has degree exactly $(d+1)/2$.
\end{definition}

\begin{lemma}
\label{lem:rsIII}
With reference to Assumption \ref{def:a3}, both
\begin{eqnarray*}
2 u b b^*+ v &=& \theta_0  \theta^*_0 + 
\theta_d  \theta^*_d, \\
-2 u b b^*+ v &=& \theta_0  \theta^*_d + 
\theta_d  \theta^*_0,
\end{eqnarray*}
where 
\begin{eqnarray*}
u &=& 1,
\label{eq:rdefIII}
\\
v &=& (2a-cd)(2a^*-c^*d)/2.
\label{eq:sdefIII}
\end{eqnarray*} 
\end{lemma}
\noindent {\it Proof:}
Use
Lemma \ref{lem:eig3}.
\hfill $\Box$ \\

\begin{theorem}
\label{thm:hatpIII}
With reference to Assumption \ref{def:a3},
\begin{eqnarray}
{\hat p}_i(\lambda) &=& p_i(u \lambda + v) \qquad \qquad (1 \leq i \leq d),
\label{eq:pihatIII}
\\
{\hat P}(\lambda) &=& P(u \lambda + v),
\label{eq:phatIII}
\end{eqnarray}
where $u,v$ are from Lemma
\ref{lem:rsIII}.
\end{theorem}
\noindent {\it Proof:}
Each side of 
(\ref{eq:pihatIII}) is a polynomial in
$\K\lbrack \lambda \rbrack$ with degree at most 1.
On each side of
(\ref{eq:pihatIII})
the polynomial has $\lambda $ coefficient $0$ (if $i$ is even)
and  $-1$ (if $i$ is odd).
Moreover at
 $\lambda=-2bb^*$ this polynomial takes the
value
$(\theta_0-\theta_i)(\theta^*_0-\theta^*_i)$.
Now (\ref{eq:pihatIII}) is true
and 
(\ref{eq:phatIII}) follows.
\hfill $\Box$ \\

\begin{proposition}
\label{prop:d4III}
With reference to Assumption \ref{def:a3} and
Definition
\ref{def:paltIII}, the polynomial
$\hat P$ is $D_4$-invariant.
\end{proposition}
\noindent {\it Proof:}
Follows from
Theorem
\ref{thm:main2}
and (\ref{eq:phatIII}).
\hfill $\Box$ \\

\begin{proposition}
\label{prop:lasttermIII}
With reference to Assumption \ref{def:a3} and
Definition
\ref{def:paltIII}, the following
 {\rm (i)}, {\rm (ii)} hold.
\begin{itemize}
\item[\rm (i)]
 ${\hat P}(2b b^*)= \zeta_d$;
\item[\rm (ii)]
 ${\hat P}(-2b b^*)= \zeta^\Downarrow_d$.
\end{itemize}
\end{proposition}
\noindent {\it Proof:}
Combine
Proposition
\ref{prop:lastterm},
Lemma \ref{lem:rsIII},
and 
(\ref{eq:phatIII}).
\hfill $\Box$ \\

\noindent We now consider the normalized Drinfel'd polynomial
for a Leonard system of type
 ${\rm III}^-$.

\begin{lemma}
\label{ex:lpIII}
{\rm \cite[Theorem~9.1]{NT:balanced}}
With reference to Assumption \ref{def:a3},
suppose $\Phi$ is a Leonard system.
 Let
$\lbrace \varphi_i\rbrace_{i=1}^d$
(resp. 
$\lbrace \phi_i\rbrace_{i=1}^d$) denote the corresponding
first (resp. second) split sequence, from
\cite[Section~3]{LS99}.
Then there exists $t \in {\overline \K}$ such that
for $1 \leq i \leq d$,
\begin{eqnarray*}
\varphi_i &=& 
\begin{cases}  
c c^*i (d-i+1) 
& \text{\rm if $i$ is even},  \\
t - 2 b b^* - 2(bc^*+cb^*)(i- \tfrac{d+1}{2})-c c^*(i-\tfrac{d+1}{2})^2
& \text{\rm if $i$ is odd},
    \end{cases} 
\\
\phi_i &=&
\begin{cases}  
c c^*i (d-i+1) 
& \text{\rm if $i$ is even},  \\
t + 2 b b^* + 2(bc^*-cb^*)(i- \tfrac{d+1}{2})-c c^*(i-\tfrac{d+1}{2})^2
& \text{\rm if $i$ is odd}.
    \end{cases} 
\end{eqnarray*}
For $i$ odd we have 
\begin{eqnarray*}
\varphi_i &=&  
c^{-1} c^{*-1} 
\bigl(\psi - b c^* - c b^* - c c^* (i-\tfrac{d+1}{2})\bigr)
\bigl(\psi +b c^* + c b^* + c c^* (i-\tfrac{d+1}{2})\bigr),
\\
\phi_i &=&  
c^{-1} c^{*-1} 
\bigl(\psi + b c^* - c b^* - c c^* (i-\tfrac{d+1}{2})\bigr)
\bigl(\psi -b c^* + c b^* + c c^* (i-\tfrac{d+1}{2})\bigr),
\end{eqnarray*}
where $\psi \in {\overline \K}$ is a solution to
\begin{eqnarray*}
\psi^2 =  t c c^* + b^2 c^{*2} + b^{* 2} c^2.
\end{eqnarray*}
\end{lemma} 

\begin{proposition}
\label{prop:rootsIII}
With reference to Assumption \ref{def:a3},
suppose $\Phi$ is a Leonard system
  and
let $\hat P$ denote the corresponding normalized Drinfel'd
polynomial.
Then the roots of
$\hat P$ are 
\begin{eqnarray*}
t + 2 \psi (i- \tfrac{d+1}{2})+ c c^* (i- \tfrac{d+1}{2})^2
\qquad \qquad (1 \leq i \leq d, \; i \;\mbox{\rm odd}),
\end{eqnarray*}
where $\psi, t $ are from Lemma \ref{ex:lpIII}.
\end{proposition}
\noindent {\it Proof:}
Abbreviate
$N=(d+1)/2$.
We begin as in the proof of
Proposition
\ref{prop:rootsI}.
We evaluate 
(\ref{eq:adj}) using
(\ref{eq:pinormd3}) 
and 
the data in Lemma 
\ref{ex:lpIII}.
The result is a hypergeometric series
\begin{eqnarray}
 {}_3F_2 \Biggl\lbrack{{-N,\;\; \;
\tfrac{x-bc^*-cb^*}{2cc^*}, \;\; \;
-\tfrac{x+bc^*+cb^*}{2cc^*}}\atop
{\tfrac{1-d}{4}-\tfrac{bc^*+cb^*-\psi}{2cc^*},\;\;
\tfrac{1-d}{4}-\tfrac{bc^*+cb^*+\psi}{2cc^*}} 
}\;  ; \; 1 \Biggr\rbrack,
\label{eq:3F2type3}
\end{eqnarray}
where $x^2= \lambda cc^*+ b^2c^{*2}+ b^{*2}c^2$.
The terms in the series
(\ref{eq:adj}) are related to the terms in
the series
(\ref{eq:3F2type3}) as follows.
For even $n$ $(0 < n < d)$ the
$n$-summand in 
(\ref{eq:adj}) plus the 
$(n-1)$-summand in 
(\ref{eq:adj}) is equal to
the $(n/2)$-summand in
(\ref{eq:3F2type3}).
The $0$-summand in
(\ref{eq:adj}) is equal to
the $0$-summand in
(\ref{eq:3F2type3}),
and the
$d$-summand in 
(\ref{eq:adj}) is equal to
the $N$-summand in
(\ref{eq:3F2type3}).
By the Saalsch\"{u}tz formula
\cite[p.~17]{GR} the sum
(\ref{eq:3F2type3}) is equal to
\begin{eqnarray}
\frac{
\bigl(\tfrac{\psi-x}{2cc^*}+\tfrac{1-d}{4}\bigr)_N
\bigl(\tfrac{\psi+x}{2cc^*}+\tfrac{1-d}{4}\bigr)_N 
}
{
\bigl(\tfrac{\psi+bc^*+cb^*}{2cc^*}+\tfrac{1-d}{4}\bigr)_N
\bigl(\tfrac{\psi-bc^*-cb^*}{2cc^*}+\tfrac{1-d}{4}\bigr)_N
}.
\label{eq:ssformIItype3}
\end{eqnarray}
The numerator in
(\ref{eq:ssformIItype3}) 
can be expressed as
\begin{eqnarray}
\label{eq:numpreIItype3}
  \prod_{\stackrel{1 \leq i \leq d}{i \;{\mbox{\rm {\tiny odd}}}}}
\tfrac{1}{4}\bigl(\tfrac{\psi-x}{cc^*}+i-\tfrac{d+1}{2}\bigr)
\bigl(\tfrac{\psi+x}{cc^*}+i-\tfrac{d+1}{2}\bigr). 
\end{eqnarray}
For odd $i$ $(1 \leq i \leq d)$ the $i$-factor in
(\ref{eq:numpreIItype3}) is 
equal to $(4cc^*)^{-1}$ times
\begin{eqnarray*}
 t +2\psi (i- \tfrac{d+1}{2})+ c c^* (i- \tfrac{d+1}{2})^2-\lambda.
\end{eqnarray*}
Therefore the numerator in (\ref{eq:ssformIItype3}) 
is a nonzero scalar multiple
of 
\begin{eqnarray*}
  \prod_{\stackrel{1 \leq i \leq d}{i \; {\mbox{\rm {\tiny odd}}}}}
\bigl(t +2\psi (i- \tfrac{d+1}{2})+ c c^* (i- \tfrac{d+1}{2})^2-\lambda \bigr).
\end{eqnarray*}
The result follows.
\hfill $\Box$ \\

\section{The normalized Drinfel'd polynomial for type
${\rm IV}$ }

In this section we discuss the
normalized Drinfel'd polynomial for a
sharp tridiagonal system of type IV.
For this type it will turn out that the normalized 
Drinfel'd polynomial is the same as the Drinfel'd polynomial,
but for notational consistency we will continue
the hat notation from Sections 9--11.

\begin{assumption}
\label{def:a4}
\rm
Let $\Phi$ denote a sharp tridiagonal system
over $\K$ that has type IV. Let
$\lbrace \theta_i\rbrace_{i=0}^3$
(resp.
$\lbrace \theta^*_i\rbrace_{i=0}^3$)
denote the eigenvalue sequence
(resp. dual eigenvalue sequence) of $\Phi$.
Let the scalars
$a,b,c,a^*,b^*,c^*$ be as in
Lemma
\ref{lem:eig4}.
\end{assumption}

\begin{definition}
\label{def:pialtIV}
\rm
With reference to Assumption \ref{def:a4},
 define ${\hat p}_i \in \K\lbrack \lambda \rbrack$ by
\begin{eqnarray*}
{\hat p}_1 &=& ab^*+ba^*+ (a+b+c)(a^*+ b^*+c^*) + \lambda,
\\
{\hat p}_2 &=& cc^*,
\\
{\hat p}_3 &=& ab^*+ba^*+ (a+b)(a^*+b^*) + \lambda.
\end{eqnarray*}
By Lemma
\ref{lem:oioiIV},
\begin{eqnarray}
{\hat p}_i(ab^*+ba^*)= (\theta_0-\theta_i)(\theta^*_0-\theta^*_i).
\label{eq:hatpnormIV}
\end{eqnarray} 
\end{definition}

\begin{definition}
\label{def:paltIV}
\rm
With reference to Assumption \ref{def:a4},
 define 
\begin{eqnarray}
\label{eq:drinnormIV}
{\hat P} = {\hat p}_1 {\hat p}_2 {\hat p}_3
+
\zeta_1 
 {\hat p}_2 {\hat p}_3
+
\zeta_2 
  {\hat p}_3
+
\zeta_3,
\end{eqnarray}
where $\lbrace \zeta_i\rbrace_{i=0}^3$ is the split
sequence of $\Phi$ and the ${\hat p}_i$ are from
Definition
\ref{def:pialtIV}. We call
$\hat P$ the {\it normalized Drinfel'd polynomial} for $\Phi$.
We note that ${\hat P}$ has degree exactly 2.
\end{definition}

\begin{lemma}
\label{thm:hatpIV}
With reference to Assumption \ref{def:a4},
\begin{eqnarray}
{\hat p}_i &=& p_i \qquad \qquad (1 \leq i \leq 3),
\label{eq:pihatIV}
\\
{\hat P} &=& P.
\label{eq:phatIV}
\end{eqnarray}
\end{lemma}
\noindent {\it Proof:}
By Definition
\ref{def:aap} and since
$\mbox{\rm Char}(\K)=2$,
\begin{eqnarray*}
p_1 &=& \theta_0 \theta^*_3 + \theta_3 \theta^*_0 +\lambda 
+(\theta_0-\theta_1)(\theta^*_0-\theta^*_1),
\\
p_2 &=&
(\theta_0-\theta_2)(\theta^*_0-\theta^*_2),
\\
p_3 &=& \theta_0 \theta^*_3 + \theta_3 \theta^*_0 +\lambda 
+(\theta_0-\theta_3)(\theta^*_0-\theta^*_3).
\end{eqnarray*}
Evaluating these lines using Lemma
\ref{lem:eig4}, Lemma
\ref{lem:oioiIV}
and comparing the result with
Definition
\ref{def:pialtIV}, we get
(\ref{eq:pihatIV}).
Line
(\ref{eq:phatIV}) follows in view of
(\ref{eq:drinp}) 
and
(\ref{eq:drinnormIV}).
\hfill $\Box$ \\

\begin{proposition}
\label{prop:d4IV}
With reference to Assumption \ref{def:a4}
and Definition
\ref{def:paltIV}, the polynomial
${\hat P}$ is 
$D_4$-invariant.
\end{proposition}
\noindent {\it Proof:}
Clear from 
Theorem
\ref{thm:main2} and
Lemma \ref{thm:hatpIV}.
\hfill $\Box$ \\

\begin{proposition}
\label{prop:lasttermIV}
With reference to Assumption \ref{def:a4}
and Definition
\ref{def:paltIV}, 
the 
following {\rm (i)}, {\rm (ii)} hold.
\begin{itemize}
\item[\rm (i)]
 ${\hat P}(a a^*+bb^*)= \zeta_3$;
\item[\rm (ii)]
 ${\hat P}(a b^*+ba^*)= \zeta^\Downarrow_3$.
\end{itemize}
\end{proposition}
\noindent {\it Proof:}
Evaluate Proposition
\ref{prop:lastterm}
using Lemma
\ref{lem:eig4}
and 
(\ref{eq:phatIV}).
\hfill $\Box$ \\

\noindent We now consider the normalized Drinfel'd polynomial
for a Leonard system of type
 ${\rm IV}$.

\begin{lemma}
\label{ex:lpIV}
{\rm \cite[Theorem~10.1]{NT:balanced}}
With reference to Assumption \ref{def:a4},
suppose $\Phi$ is a Leonard system.
 Let
$\lbrace \varphi_i\rbrace_{i=1}^3$
(resp. 
$\lbrace \phi_i\rbrace_{i=1}^3$) denote the corresponding
first (resp. second) split sequence, from
\cite[Section~3]{LS99}. 
Then there exists $\varphi \in \K$ such that
\begin{eqnarray*}
&&\varphi_1 = \varphi, \qquad 
\varphi_2 = c c^*, \qquad
\varphi_3 = \varphi + (a+b) c^*+ c(a^*+b^*),
\\
&&\phi_1 = \varphi + (a+b)(a^*+b^*+c^*), \qquad 
\phi_2 = c c^*, \qquad
\phi_3 = \varphi + 
(a+b+c)(a^*+b^*).
\end{eqnarray*} 
\end{lemma} 

\begin{proposition}
\label{prop:polyIV}
With reference to Assumption \ref{def:a4},
suppose $\Phi$ is a Leonard system
  and
let $\hat P$ denote the corresponding normalized Drinfel'd
polynomial.
Then
$\hat P$ is $cc^*$ times 
\begin{eqnarray*}
&&(\lambda + a b^*+b a^*)^2
\; + \; 
(\lambda + a b^*+b a^*)(ac^* + bc^* + ca^* + cb^* + c c^*)
\\
&& \qquad \qquad 
+ \;
\varphi^2 \; + \;  \varphi (ac^* + bc^* + c a^*+ cb^*)\;  + \;
 (a+b)(a^*+b^*)(a+b+c)(a^*+b^*+c^*).
\end{eqnarray*}
\end{proposition}
\noindent {\it Proof:}
Evaluate
(\ref{eq:drinnormIV}) using
Definition
\ref{def:pialtIV}, Lemma 
\ref{ex:lpIV},
 and
$\zeta_1=\varphi_1$,
$\zeta_2=\varphi_1 \varphi_2$,
$\zeta_3=\varphi_1 \varphi_2 \varphi_3$.
\hfill $\Box$ \\

\noindent Referring to 
Proposition
\ref{prop:polyIV},
we caution the reader that since 
$\mbox{\rm Char}(\K)=2$ the roots of $\hat P$ cannot
be obtained using the quadratic formula. To get these
roots we proceed as follows.

\begin{lemma}
\label{ex:lpIVpsi}
With reference to Assumption \ref{def:a4},
suppose $\Phi$ is a Leonard system.
 Let
$\lbrace \varphi_i\rbrace_{i=1}^3$
(resp. 
$\lbrace \phi_i\rbrace_{i=1}^3$) denote the corresponding
first (resp. second) split sequence. 
Then there exists $\psi \in {\overline \K}$ such that
\begin{eqnarray*}
c c^*\varphi_1 &=& (a c^*+a^* c+c c^* \psi)
 (b c^*+b^* c+c c^* + c c^* \psi),
\\
c c^*\varphi_3 &=& (b c^*+b^* c+c c^* \psi)
 (a c^*+a^* c+c c^* + c c^* \psi),
\\
c c^*\phi_1 &=& (b c^*+a^* c+c c^* \psi)
 (a c^*+b^* c+c c^* + c c^* \psi),
\\
c c^*\phi_3 &=& (a c^*+b^* c+c c^* \psi)
 (b c^*+a^* c+c c^* + c c^* \psi).
\end{eqnarray*}
\end{lemma}
\noindent {\it Proof:}
Since $cc^*\not=0$
and since $\overline \K$ is algebraically closed
there exists $\psi \in {\overline \K}$
that satisfies the first of the four equations above.
The remaining three equations follow
in view of the data in Lemma
\ref{ex:lpIV}.
\hfill $\Box$ \\

\noindent We comment on the uniqueness of the
scalar $\psi$ in Lemma
\ref{ex:lpIVpsi}.

\begin{note}
\label{ex:lpIVpsiunique}
\rm
If $\psi \in {\overline \K}$ satisfies
the four equations in Lemma
\ref{ex:lpIVpsi},
then
\begin{eqnarray*}
\psi + \frac{a+b}{c}+\frac{a^*+b^*}{c^*}+1
\end{eqnarray*}
satisfies these equations 
and no other scalar in $\overline \K$ satisfies these equations. 
\end{note}

\begin{theorem}
\label{ex:lpIVpsiroots}
With reference to Assumption \ref{def:a4},
suppose $\Phi$ is a Leonard system
  and
let $\hat P$ denote the corresponding normalized Drinfel'd
polynomial.
Then the roots of $\hat P$ are
\begin{eqnarray*}
&&
a b^* + b a^* + (c c^* \psi + a c^* +  b^*c)(c c^* \psi + a^* c + b c^*)c^{-1}
c^{*-1},
\\
&&
a b^* + b a^* + (c c^* \psi + a c^* +  b^*c + c c^*)(c c^* \psi + a^* c + b c^*
+ c c^*)c^{-1}
c^{*-1},
\end{eqnarray*}
where $\psi$ is from
Lemma
\ref{ex:lpIVpsi}.
\end{theorem}
\noindent {\it Proof:}
Denote the above expressions by $z_1$, $z_2$.
One finds ${\hat P}(z_1)=0$,
 ${\hat P}(z_2)=0$ using
Proposition
\ref{prop:polyIV}
and Lemma
\ref{ex:lpIVpsi}.
\hfill $\Box$ \\

\noindent 
Combining 
Propositions
\ref{prop:d4I},
\ref{prop:d4II},
\ref{prop:d4III},
\ref{prop:d4IV}
and the second paragraph in
Section 11,
we obtain the following theorem.

\begin{theorem}
\label{thm:maind4n}
Let $\Phi$ denote a sharp tridiagonal system
and let $\hat P$ denote the corresponding normalized
 Drinfel'd polynomial.
Then $\hat P$ is $D_4$-invariant.
\end{theorem}

\begin{definition}
\label{def:normdrinpoly}
\rm
Let $A,A^*$ denote a sharp tridiagonal pair.
By the {\it normalized Drinfel'd polynomial}
for $A,A^*$ we mean the normalized Drinfel'd polynomial
for an associated tridiagonal system.
By construction
$A,A^*$ and $A^*,A$ have the same normalized Drinfel'd
 polynomial.
\end{definition}

\section{Why $P$ is called the Drinfel'd polynomial}

\noindent 
Let $A,A^*$ denote a sharp tridiagonal pair.
Earlier in the paper, we associated with this pair 
a polynomial $P$ called the Drinfel'd polynomial.
In this section
we justify the name
by relating $P$ to the classical
Drinfel'd polynomial from the theory of
Lie algebras and quantum groups.

\medskip
\noindent 
Throughout this section we assume that 
the field $\K$ is 
algebraically closed with characteristic zero.

\begin{definition}
\label{def:krawtdp}
\rm
\cite[Section~1] {IT:Krawt}
Let $A,A^*$ denote a 
tridiagonal pair over $\F$ that has diameter $d$. This pair
is said to have {\it Krawtchouk type}
whenever the sequence 
$\lbrace d-2i\rbrace_{i=0}^d$ is a standard ordering of the eigenvalues
of $A$ and a standard ordering of the eigenvalues of $A^*$.
\end{definition}

\noindent Let $A,A^*$ denote a tridiagonal pair on $V$ that has
Krawtchouk type. 
By 
\cite[Theorem 1.8]{Ha} 
the pair $A,A^*$ induces on $V$
a module structure for the 3-point 
$\mathfrak{sl}_2$ loop algebra
\cite[Definition 1.1]{HT}.
Associated with this module is a Drinfel'd polynomial
\cite[Definition 9.13]{Ev},
\cite[Lemma 13.2]{IT:Krawt} 
which we denote by $P_{A,A^*}$. In the notation of
the present paper 
$P_{A,A^*}$ looks as follows.

\begin{definition}
\label{def:drinkrlit}
\rm
\cite[Definition~13.1]{IT:Krawt}
Let $A,A^*$ denote a tridiagonal pair over $\K$
that has Krawtchouk type. Define $P_{A,A^*} \in \K\lbrack \lambda \rbrack$
by 
\begin{eqnarray}
\label{eq:drinkrlit}
P_{A,A^*} = \sum_{i=0}^d \frac{(-1)^i \zeta_i \lambda^i}{(i!)^2 4^i},
\end{eqnarray}
where $\lbrace \zeta_i \rbrace_{i=0}^d $ is the split sequence
of $A,A^*$ associated with  
the standard ordering $\lbrace d-2i\rbrace_{i=0}^d$
(resp.  
 $\lbrace 2i-d\rbrace_{i=0}^d$)
of the eigenvalues of $A$ (resp. $A^*$).
\end{definition}

\begin{theorem}
\label{thm:krsame}
Let $A,A^*$ denote a tridiagonal pair over $\K$ that has Krawtchouk
type, and let $\hat P$ denote the associated normalized Drinfel'd 
polynomial from Definition
\ref{def:normdrinpoly}.
Then
\begin{eqnarray}
{\hat P}(\lambda) = 
(-1)^d (d !)^2 (\lambda + 2)^d P_{A,A^*}(4(\lambda+2)^{-1}),
\label{eq:krsame}
\end{eqnarray}
where $P_{A,A^*}$ is from Definition
\ref{def:drinkrlit}
 and $d$ is the diameter of $A,A^*$.
\end{theorem}
\noindent {\it Proof:}
To describe ${\hat P}$ we associate 
with $A,A^*$ a tridiagonal system.
 For $0 \leq i \leq d$
let $E_i$ (resp. $E^*_i$)
denote the primitive idempotent of $A$ (resp. $A^*$)
associated with the eigenvalue $d-2i$
(resp. $2i-d$). Then
 $\Phi=(A;\lbrace E_i\rbrace_{i=0}^d;
A^*; \lbrace E^*_i\rbrace_{i=0}^d)$ 
is a tridiagonal system with eigenvalue sequence
$\lbrace d-2i\rbrace_{i=0}^d$
and dual eigenvalue sequence
$\lbrace 2i-d\rbrace_{i=0}^d$.
By Definition 
     \ref{defe}
$\Phi$ is type $\rm II$, and the equations of
Lemma
\ref{lem:eig2} are satisfied by
\begin{eqnarray*}
a=0, \quad 
a^*=0, \quad 
b=-2, \quad 
b^*=2, \quad 
c=0, \quad 
c^*=0.
\end{eqnarray*}
Evaluating 
(\ref{eq:normdp2}) using
this we find ${\hat p}_i(\lambda) = -i^2(\lambda+2)$
for $1 \leq i \leq d$. Now using
(\ref{eq:drinnormII}),
\begin{eqnarray}
\label{eq:htpkr}
{\hat P}(\lambda) = (-1)^d (d !)^2 (\lambda+2)^d \sum_{i=0}^d 
\frac{(-1)^i \zeta_i}{(i!)^2 (\lambda+2)^i}.
\end{eqnarray}
Comparing 
(\ref{eq:drinkrlit})
and
(\ref{eq:htpkr}) we obtain
(\ref{eq:krsame}).
\hfill $\Box$ \\

\noindent
For the rest of this section
fix a nonzero $q \in \K$
that is not a root of $1$. For all integers $n\geq 0$ define
\begin{eqnarray*}
\lbrack n \rbrack_q &=& \frac{q^n-q^{-n}}{q-q^{-1}},
\\
\lbrack n \rbrack^!_q &=& 
\lbrack n \rbrack_q 
\lbrack n-1 \rbrack_q 
\cdots 
\lbrack 1 \rbrack_q.
\end{eqnarray*}
We interpret 
$\lbrack 0 \rbrack^!_q = 1$.

\begin{definition}
\label{def:geom}
\rm
{\rm
 \cite[Definition 2.6]{NN}}
Let $A,A^*$ denote a tridiagonal pair over $\K$ that has diameter $d$.
Then $A,A^*$ 
is called {\it $q$-geometric} whenever the sequence
$\lbrace q^{d-2i}\rbrace_{i=0}^d$ is a standard ordering of the eigenvalues
of $A$ and a standard ordering of the eigenvalues of $A^*$.
\end{definition}

\noindent Let $A,A^*$ denote a $q$-geometric 
tridiagonal pair on $V$.
By \cite[Theorem~3.3,~Theorem~13.1]{tdanduq}
the pair $A,A^*$ induces on $V$
a module structure for 
the quantum group 
$U_q({\widehat{\mathfrak{sl}}}_2)$.
Associated with this module is a Drinfel'd polynomial
 \cite[Definition 4.2]{NN}
which we will denote by $P_{A,A^*}$. In the notation of
the present paper 
$P_{A,A^*}$ looks as follows.

\begin{definition}
\label{def:drinnnlit}
\rm
\cite[Definition 4.2]{NN}
Let $A,A^*$ denote a $q$-geometric tridiagonal pair over $\K$.
Define $P_{A,A^*} \in \K\lbrack \lambda \rbrack$
by 
\begin{eqnarray}
\label{eq:dringeomlit}
P_{A,A^*} = \sum_{i=0}^d \frac{(-1)^i \zeta_i q^i\lambda^i}
{(\lbrack i \rbrack^!_q)^2},
\end{eqnarray}
where $\lbrace \zeta_i \rbrace_{i=0}^d $ is the split sequence
of $A,A^*$ associated with  
the standard ordering $\lbrace q^{2i-d}\rbrace_{i=0}^d$
(resp.  
 $\lbrace q^{d-2i}\rbrace_{i=0}^d$)
of the eigenvalues of $A$ (resp. $A^*$).
\end{definition}

\begin{theorem}
\label{thm:geomsame}
Let $A,A^*$ denote a $q$-geometric tridiagonal pair over $\K$,
 and let $\hat P$ denote the associated normalized Drinfel'd 
polynomial from Definition
\ref{def:normdrinpoly}.
Then
\begin{eqnarray}
{\hat P}(\lambda) = 
(-1)^d (\lbrack d\rbrack^!_q)^2 (q-q^{-1})^{2d} \lambda^d
 P_{A,A^*}(\lambda^{-1}q^{-1}(q-q^{-1})^{-2}),
\label{eq:geomsame}
\end{eqnarray}
where $P_{A,A^*}$ is from
 Definition
\ref{def:drinnnlit}
 and $d$ is the diameter of $A,A^*$.
\end{theorem}
\noindent {\it Proof:}
To describe ${\hat P}$ we associate with
 $A,A^*$ 
a tridiagonal system.
 For $0 \leq i \leq d$
let $E_i$ (resp. $E^*_i$)
denote the primitive idempotent of $A$ (resp. $A^*$)
associated with the eigenvalue $q^{2i-d}$
(resp. $q^{d-2i}$). Then
 $\Phi=(A;\lbrace E_i\rbrace_{i=0}^d;
A^*; \lbrace E^*_i\rbrace_{i=0}^d)$ 
is a tridiagonal system with eigenvalue sequence
$\lbrace q^{2i-d}\rbrace_{i=0}^d$
and dual eigenvalue sequence
$\lbrace q^{d-2i}\rbrace_{i=0}^d$.
Referring to 
Definition
\ref{def:base},
 we have $\beta=q^2+q^{-2}$ and $q$ is not
a root of unity so  
$\Phi$ is type $\rm I$.
Without loss we may take the scalar $q$ from  
Lemma
\ref{lem:eig1} to be the present $q$ that we fixed
above Definition
\ref{def:geom}. The equations 
of Lemma \ref{lem:eig1}
are satisfied by
\begin{eqnarray*}
a=0, \quad 
a^*=0, \quad 
b=1, \quad 
b^*=0, \quad 
c=0, \quad 
c^*=1.
\end{eqnarray*}
Evaluating 
(\ref{eq:pinormcase1})
using this we find
${\hat p}_i=-(q^i-q^{-i})^2\lambda $
for $1 \leq i \leq d$. 
 Now using
(\ref{eq:drinnorm}),
\begin{eqnarray}
\label{eq:htpgeom}
{\hat P}(\lambda) = (-1)^d
(\lbrack d \rbrack^!_q)^2(q-q^{-1})^{2d} \lambda ^d
\sum_{i=0}^d 
\frac{(-1)^i \zeta_i}{(q-q^{-1})^{2i}(\lbrack i \rbrack^!_q)^2\lambda^i}.
\end{eqnarray}
Comparing 
(\ref{eq:dringeomlit})
and
(\ref{eq:htpgeom}) we  obtain
(\ref{eq:geomsame}).
\hfill $\Box$ \\

\section{Acknowledgements}

\indent
The authors thank Brian Curtin, Eric Egge, and Kazumasa Nomura
for giving
this paper a close reading and offering many valuable suggestions.

\noindent 
Tatsuro Ito
 \hfil\break
\noindent 
Division of Mathematical and Physical Sciences 
 \hfil\break
\noindent 
Graduate School of Natural Science and Technology
 \hfil\break
\noindent 
Kanazawa University 
 \hfil\break
\noindent 
Kakuma-machi, Kanazawa 920-1192, Japan 
 \hfil\break
\noindent 
email: {\tt tatsuro@kenroku.kanazawa-u.ac.jp}

\bigskip

\noindent Paul Terwilliger \hfil\break
\noindent Department of Mathematics \hfil\break
\noindent University of Wisconsin \hfil\break
\noindent 480 Lincoln Drive \hfil\break
\noindent Madison, WI 53706-1388 USA \hfil\break
\noindent email: {\tt terwilli@math.wisc.edu }\hfil\break

\end{document}